

\input amstex
\expandafter\ifx\csname mathdefs.tex\endcsname\relax
  \expandafter\gdef\csname mathdefs.tex\endcsname{}
\else \message{Hey!  Apparently you were trying to
  \string twice.   This does not make sense.} 
\errmessage{Please edit your file (probably \jobname.tex) and remove
any duplicate ``\string\input'' lines} \fi




\catcode`\X=12\catcode`\@=11

\def\n@wcount{\alloc@0\count\countdef\insc@unt}
\def\n@wwrite{\alloc@7\write\chardef\sixt@@n}
\def\n@wread{\alloc@6\read\chardef\sixt@@n}
\def\r@s@t{\relax}\def\v@idline{\par}\def\@mputate#1/{#1}
\def\l@c@l#1X{\firstpart.#1}\def\gl@b@l#1X{#1}\def\t@d@l#1X{{}}

\def\crossrefs#1{\ifx\all#1\let\tr@ce=\all\else\def\tr@ce{#1,}\fi
   \n@wwrite\cit@tionsout\openout\cit@tionsout=\jobname.cit 
   \write\cit@tionsout{\tr@ce}\expandafter\setfl@gs\tr@ce,}
\def\setfl@gs#1,{\def\@{#1}\ifx\@\empty\let\next=\relax
   \else\let\next=\setfl@gs\expandafter\xdef
   \csname#1tr@cetrue\endcsname{}\fi\next}
\def\m@ketag#1#2{\expandafter\n@wcount\csname#2tagno\endcsname
     \csname#2tagno\endcsname=0\let\tail=\all\xdef\all{\tail#2,}
   \ifx#1\l@c@l\let\tail=\r@s@t\xdef\r@s@t{\csname#2tagno\endcsname=0\tail}\fi
   \expandafter\gdef\csname#2cite\endcsname##1{\expandafter
     \ifx\csname#2tag##1\endcsname\relax?\else\csname#2tag##1\endcsname\fi
     \expandafter\ifx\csname#2tr@cetrue\endcsname\relax\else
     \write\cit@tionsout{#2tag ##1 cited on page \folio.}\fi}
   \expandafter\gdef\csname#2page\endcsname##1{\expandafter
     \ifx\csname#2page##1\endcsname\relax?\else\csname#2page##1\endcsname\fi
     \expandafter\ifx\csname#2tr@cetrue\endcsname\relax\else
     \write\cit@tionsout{#2tag ##1 cited on page \folio.}\fi}
   \expandafter\gdef\csname#2tag\endcsname##1{\expandafter
      \ifx\csname#2check##1\endcsname\relax
      \expandafter\xdef\csname#2check##1\endcsname{}%
      \else\immediate\write16{Warning: #2tag ##1 used more than once.}\fi
      \multit@g{#1}{#2}##1/X%
      \write\t@gsout{#2tag ##1 assigned number \csname#2tag##1\endcsname\space
      on page \number\count0.}%
   \csname#2tag##1\endcsname}}
\def\multit@g#1#2#3/#4X{\def\t@mp{#4}\ifx\t@mp\empty%
      \global\advance\csname#2tagno\endcsname by 1 
      \expandafter\xdef\csname#2tag#3\endcsname
      {#1\number\csname#2tagno\endcsnameX}%
   \else\expandafter\ifx\csname#2last#3\endcsname\relax
      \expandafter\n@wcount\csname#2last#3\endcsname
      \global\advance\csname#2tagno\endcsname by 1 
      \expandafter\xdef\csname#2tag#3\endcsname
      {#1\number\csname#2tagno\endcsnameX}
      \write\t@gsout{#2tag #3 assigned number \csname#2tag#3\endcsname\space
      on page \number\count0.}\fi
   \global\advance\csname#2last#3\endcsname by 1
   \def\t@mp{\expandafter\xdef\csname#2tag#3/}%
   \expandafter\t@mp\@mputate#4\endcsname
   {\csname#2tag#3\endcsname\lastpart{\csname#2last#3\endcsname}}\fi}
\def\t@gs#1{\def\all{}\m@ketag#1e\m@ketag#1s\m@ketag\t@d@l p
   \m@ketag\gl@b@l r \n@wread\t@gsin
   \openin\t@gsin=\jobname.tgs \re@der \closein\t@gsin
   \n@wwrite\t@gsout\openout\t@gsout=\jobname.tgs }
\outer\def\localtags{\t@gs\l@c@l}
\outer\def\globaltags{\t@gs\gl@b@l}
\outer\def\newlocaltag#1{\m@ketag\l@c@l{#1}}
\outer\def\newglobaltag#1{\m@ketag\gl@b@l{#1}}

\newif\ifpr@ 
\def\m@kecs #1tag #2 assigned number #3 on page #4.%
   {\expandafter\gdef\csname#1tag#2\endcsname{#3}
   \expandafter\gdef\csname#1page#2\endcsname{#4}
   \ifpr@\expandafter\xdef\csname#1check#2\endcsname{}\fi}
\def\re@der{\ifeof\t@gsin\let\next=\relax\else
   \read\t@gsin to\t@gline\ifx\t@gline\v@idline\else
   \expandafter\m@kecs \t@gline\fi\let \next=\re@der\fi\next}
\def\pretags#1{\pr@true\pret@gs#1,,}
\def\pret@gs#1,{\def\@{#1}\ifx\@\empty\let\n@xtfile=\relax
   \else\let\n@xtfile=\pret@gs \openin\t@gsin=#1.tgs \message{#1} \re@der 
   \closein\t@gsin\fi \n@xtfile}

\newcount\sectno\sectno=0\newcount\subsectno\subsectno=0
\newif\ifultr@local \def\ultralocal{\ultr@localtrue}
\def\firstpart{\number\sectno}
\def\lastpart#1{\ifcase#1 \or a\or b\or c\or d\or e\or f\or g\or h\or 
   i\or k\or l\or m\or n\or o\or p\or q\or r\or s\or t\or u\or v\or w\or 
   x\or y\or z \fi}

\def\resetall{\global\advance\sectno by 1\subsectno=0
   \gdef\firstpart{\number\sectno}\r@s@t}
\def\resetsub{\global\advance\subsectno by 1
   \gdef\firstpart{\number\sectno.\number\subsectno}\r@s@t}
\def\newsection#1\par{\resetall\vskip0pt plus.3\vsize\penalty-250
   \vskip0pt plus-.3\vsize\bigskip\bigskip
   \message{#1}\leftline{\bf#1}\nobreak\bigskip}
\def\subsection#1\par{\ifultr@local\resetsub\fi
   \vskip0pt plus.2\vsize\penalty-250\vskip0pt plus-.2\vsize
   \bigskip\smallskip\message{#1}\leftline{\bf#1}\nobreak\medskip}

\def\t@gsoff#1,{\def\@{#1}\ifx\@\empty\let\next=\relax\else\let\next=\t@gsoff
   \def\@@{p}\ifx\@\@@\else
   \expandafter\gdef\csname#1cite\endcsname##1{\zeigen{##1}}
   \expandafter\gdef\csname#1page\endcsname##1{?}
   \expandafter\gdef\csname#1tag\endcsname##1{\zeigen{##1}}\fi\fi\next}
\def\verbatimtags{\ifx\all\relax\else\expandafter\t@gsoff\all,\fi}
\def\zeigen#1{\hbox{$\langle$}#1\hbox{$\rangle$}}

\def\(#1){\edef\dot@g{\ifmmode\ifinner(\hbox{\noexpand\etag{#1}})
   \else\noexpand\eqno(\hbox{\noexpand\etag{#1}})\fi
   \else(\noexpand\ecite{#1})\fi}\dot@g}

\newif\ifbr@ck
\def\eat#1{}
\def\[#1]{\br@cktrue[\br@cket#1'X]}
\def\br@cket#1'#2X{\def\temp{#2}\ifx\temp\empty\let\next\eat
   \else\let\next\br@cket\fi
   \ifbr@ck\br@ckfalse\br@ck@t#1,X\else\br@cktrue#1\fi\next#2X}
\def\br@ck@t#1,#2X{\def\temp{#2}\ifx\temp\empty\let\neext\eat
   \else\let\neext\br@ck@t\def\temp{,}\fi
   \def\teemp{#1}\ifx\teemp\empty\else\rcite{#1}\fi\temp\neext#2X}
\def\resetbr@cket{\gdef\[##1]{[\rtag{##1}]}}
\def\references{\resetbr@cket\newsection References\par}

\newtoks\symb@ls\newtoks\s@mb@ls\newtoks\p@gelist\n@wcount\ftn@mber
    \ftn@mber=1\newif\ifftn@mbers\ftn@mbersfalse\newif\ifbyp@ge\byp@gefalse
\def\defm@rk{\ifftn@mbers\n@mberm@rk\else\symb@lm@rk\fi}
\def\n@mberm@rk{\xdef\m@rk{{\the\ftn@mber}}%
    \global\advance\ftn@mber by 1 }
\def\rot@te#1{\let\temp=#1\global#1=\expandafter\r@t@te\the\temp,X}
\def\r@t@te#1,#2X{{#2#1}\xdef\m@rk{{#1}}}
\def\b@@st#1{{$^{#1}$}}\def\str@p#1{#1}
\def\symb@lm@rk{\ifbyp@ge\rot@te\p@gelist\ifnum\expandafter\str@p\m@rk=1 
    \s@mb@ls=\symb@ls\fi\write\f@nsout{\number\count0}\fi \rot@te\s@mb@ls}
\def\byp@ge{\byp@getrue\n@wwrite\f@nsin\openin\f@nsin=\jobname.fns 
    \n@wcount\currentp@ge\currentp@ge=0\p@gelist={0}
    \re@dfns\closein\f@nsin\rot@te\p@gelist
    \n@wread\f@nsout\openout\f@nsout=\jobname.fns }
\def\m@kelist#1X#2{{#1,#2}}
\def\re@dfns{\ifeof\f@nsin\let\next=\relax\else\read\f@nsin to \f@nline
    \ifx\f@nline\v@idline\else\let\t@mplist=\p@gelist
    \ifnum\currentp@ge=\f@nline
    \global\p@gelist=\expandafter\m@kelist\the\t@mplistX0
    \else\currentp@ge=\f@nline
    \global\p@gelist=\expandafter\m@kelist\the\t@mplistX1\fi\fi
    \let\next=\re@dfns\fi\next}
\def\symbols#1{\symb@ls={#1}\s@mb@ls=\symb@ls} 
\def\bigsymbol{\textstyle}
\symbols{\bigsymbol\ast,\dagger,\ddagger,\sharp,\flat,\natural,\star}
\def\ftnumbers{\ftn@mberstrue} \def\ftsymbols{\ftn@mbersfalse}
\def\paginal{\byp@ge} \def\resetftnumbers{\ftn@mber=1}
\def\ftnote#1{\defm@rk\expandafter\expandafter\expandafter\footnote
    \expandafter\b@@st\m@rk{#1}}

\long\def\jump#1\endjump{}
\def\ssum{\mathop{\lower .1em\hbox{$\textstyle\Sigma$}}\nolimits}

\def\qed{\nobreak\kern 1em \vrule height .5em width .5em depth 0em}
\def\newneq{\hbox{\rlap{\hbox to 1\wd9{\hss$=$\hss}}\raise .1em 
   \hbox to 1\wd9{\hss$\scriptscriptstyle/$\hss}}}
\def\subsetne{\setbox9 = \hbox{$\subset$}\mathrel{\hbox{\rlap
   {\lower .4em \newneq}\raise .13em \hbox{$\subset$}}}}
\def\supsetne{\setbox9 = \hbox{$\subset$}\mathrel{\hbox{\rlap
   {\lower .4em \newneq}\raise .13em \hbox{$\supset$}}}}

\def\vbar{\mathchoice{\vrule height6.3ptdepth-.5ptwidth.8pt\kern-.8pt}
   {\vrule height6.3ptdepth-.5ptwidth.8pt\kern-.8pt}
   {\vrule height4.1ptdepth-.35ptwidth.6pt\kern-.6pt}
   {\vrule height3.1ptdepth-.25ptwidth.5pt\kern-.5pt}}
\def\f@dge{\mathchoice{}{}{\mkern.5mu}{\mkern.8mu}}
\def\b@c#1#2{{\rm \mkern#2mu\vbar\mkern-#2mu#1}}
\def\b@b#1{{\rm I\mkern-3.5mu #1}}
\def\b@a#1#2{{\rm #1\mkern-#2mu\f@dge #1}}
\def\bb#1{{\count4=`#1 \advance\count4by-64 \ifcase\count4\or\b@a A{11.5}\or
   \b@b B\or\b@c C{5}\or\b@b D\or\b@b E\or\b@b F \or\b@c G{5}\or\b@b H\or
   \b@b I\or\b@c J{3}\or\b@b K\or\b@b L \or\b@b M\or\b@b N\or\b@c O{5} \or
   \b@b P\or\b@c Q{5}\or\b@b R\or\b@a S{8}\or\b@a T{10.5}\or\b@c U{5}\or
   \b@a V{12}\or\b@a W{16.5}\or\b@a X{11}\or\b@a Y{11.7}\or\b@a Z{7.5}\fi}}

\catcode`\X=11 \catcode`\@=12

\expandafter\ifx\csname citeadd.tex\endcsname\relax
\expandafter\gdef\csname citeadd.tex\endcsname{}
\else \message{Hey!  Apparently you were trying to
\string twice.   This does not make sense.} 
\errmessage{Please edit your file (probably \jobname.tex) and remove
any duplicate ``\string\input'' lines} \fi

\sectno=-2   
\localtags
\ifx\shlhetal\undefinedcontrolsequenc\let\shlhetal\relax\fi
\NoBlackBoxes
\define\mr{\medskip\roster}
\define\sn{\smallskip\noindent}
\define\mn{\medskip\noindent}
\define\bn{\bigskip\noindent}
\define\ub{\underbar}
\define\wilog{\text{without loss of generality}}
\define\ermn{\endroster\medskip\noindent}
\define\dbca{\dsize\bigcap}
\define\dbcu{\dsize\bigcup}
\define \nl{\newline}
\documentstyle {amsppt}
\topmatter
\title{Strong dichotomy of cardinality \\
Sh664} \endtitle
\rightheadtext{Strong dichotomy of cardinality}
\author {Saharon Shelah \thanks {\null\newline Research supported by the
German-Israeli Foundation for Scientific Research \null\newline
I would like to thank Alice Leonhardt for the beautiful typing. \null\newline
 Written Spring '97 \null\newline
 First Typed - 97/Sept/2 \null\newline
 Latest Revision - 98/July/28} \endthanks} \endauthor  
\affil{Institute of Mathematics\\
 The Hebrew University\\
 Jerusalem, Israel
 \medskip
 Rutgers University\\
 Mathematics Department\\
 New Brunswick, NJ  USA} \endaffil
\endtopmatter
\document  

\expandafter\ifx\csname alice2jlem.tex\endcsname\relax
  \expandafter\gdef\csname alice2jlem.tex\endcsname{}
\else \message{Hey!  Apparently you were trying to
\string  twice.   This does not make sense.}
\errmessage{Please edit your file (probably \jobname.tex) and remove
any duplicate ``\string\input'' lines} \fi

\expandafter\ifx\csname bib4plain.tex\endcsname\relax
  \expandafter\gdef\csname bib4plain.tex\endcsname{}
\else \message{Hey!  Apparently you were trying to \string twice.   This does not make sense.}
\errmessage{Please edit your file (probably \jobname.tex) and remove
any duplicate ``\string\input'' lines} \fi

\def\renewcommand{\newcommand}	       
\edef\cite{\the\catcode`@}%
\catcode`@ = 11
\let\@oldatcatcode = \cite
\chardef\@letter = 11
\chardef\@other = 12
%
%
%
%
\def\@innerdef#1#2{\edef#1{\expandafter\noexpand\csname #2\endcsname}}%
%
%
\@innerdef\@innernewcount{newcount}%
\@innerdef\@innernewdimen{newdimen}%
\@innerdef\@innernewif{newif}%
\@innerdef\@innernewwrite{newwrite}%
%
%
%
\def\@gobble#1{}%
%
%
%
\ifx\inputlineno\@undefined
   \let\@linenumber = \empty 
\else
   \def\@linenumber{\the\inputlineno:\space}%
\fi
%
%
%
\def\@futurenonspacelet#1{\def\cs{#1}%
   \afterassignment\@stepone\let\@nexttoken=
}%
\begingroup 
\def\\{\global\let\@stoken= }%
\\ 
\endgroup
\def\@stepone{\expandafter\futurelet\cs\@steptwo}%
\def\@steptwo{\expandafter\ifx\cs\@stoken\let\@@next=\@stepthree
   \else\let\@@next=\@nexttoken\fi \@@next}%
\def\@stepthree{\afterassignment\@stepone\let\@@next= }%
%
%
%
\def\@getoptionalarg#1{%
   \let\@optionaltemp = #1%
   \let\@optionalnext = \relax
   \@futurenonspacelet\@optionalnext\@bracketcheck
}%
%
%
\def\@bracketcheck{%
   \ifx [\@optionalnext
      \expandafter\@@getoptionalarg
   \else
      \let\@optionalarg = \empty
      \expandafter\@optionaltemp
   \fi
}%
\def\@@getoptionalarg[#1]{%
   \def\@optionalarg{#1}%
   \@optionaltemp
}%
%
%
%
\def\@nnil{\@nil}%
\def\@fornoop#1\@@#2#3{}%
\def\@for#1:=#2\do#3{%
   \edef\@fortmp{#2}%
   \ifx\@fortmp\empty \else
      \expandafter\@forloop#2,\@nil,\@nil\@@#1{#3}%
   \fi
}%
\def\@forloop#1,#2,#3\@@#4#5{\def#4{#1}\ifx #4\@nnil \else
       #5\def#4{#2}\ifx #4\@nnil \else#5\@iforloop #3\@@#4{#5}\fi\fi
}%
\def\@iforloop#1,#2\@@#3#4{\def#3{#1}\ifx #3\@nnil
       \let\@nextwhile=\@fornoop \else
      #4\relax\let\@nextwhile=\@iforloop\fi\@nextwhile#2\@@#3{#4}%
}%
%
%
%
\@innernewif\if@fileexists
\def\@testfileexistence{\@getoptionalarg\@finishtestfileexistence}%
\def\@finishtestfileexistence#1{%
   \begingroup
      \def\extension{#1}%
      \immediate\openin0 =
         \ifx\@optionalarg\empty\jobname\else\@optionalarg\fi
         \ifx\extension\empty \else .#1\fi
         \space
      \ifeof 0
         \global\@fileexistsfalse
      \else
         \global\@fileexiststrue
      \fi
      \immediate\closein0
   \endgroup
}%
%
%
%
%
\def\bibliographystyle#1{%
   \@readauxfile
   \@writeaux{\string\bibstyle{#1}}%
}%
\let\bibstyle = \@gobble
%
%
\let\bblfilebasename = \jobname
\def\bibliography#1{%
   \@readauxfile
   \@writeaux{\string\bibdata{#1}}%
   \@testfileexistence[\bblfilebasename]{bbl}%
   \if@fileexists
      \nobreak
      \@readbblfile
   \fi
}%
\let\bibdata = \@gobble
%
%
\def\nocite#1{%
   \@readauxfile
   \@writeaux{\string\citation{#1}}%
}%
\@innernewif\if@notfirstcitation
%
%
\def\cite{\@getoptionalarg\@cite}%
%
%
\def\@cite#1{%
   \let\@citenotetext = \@optionalarg
   \printcitestart
   \nocite{#1}%
   \@notfirstcitationfalse
   \@for \@citation :=#1\do
   {%
      \expandafter\@onecitation\@citation\@@
   }%
   \ifx\empty\@citenotetext\else
      \printcitenote{\@citenotetext}%
   \fi
   \printcitefinish
}%
\def\@onecitation#1\@@{%
   \if@notfirstcitation
      \printbetweencitations
   \fi
   \expandafter \ifx \csname\@citelabel{#1}\endcsname \relax
      \if@citewarning
         \message{\@linenumber Undefined citation `#1'.}%
      \fi
      \expandafter\gdef\csname\@citelabel{#1}\endcsname{%
\strut
\vadjust{\vskip-\dp\strutbox
\vbox to 0pt{\vss\parindent0cm \leftskip=\hsize 
\advance\leftskip3mm
\advance\hsize 4cm\strut\openup-4pt 
\rightskip 0cm plus 1cm minus 0.5cm ?  #1 ?\strut}}
         {\tt
            \escapechar = -1
            \nobreak\hskip0pt
            \expandafter\string\csname#1\endcsname
            \nobreak\hskip0pt
         }%
      }%
   \fi
   \csname\@citelabel{#1}\endcsname
   \@notfirstcitationtrue
}%
%
%
\def\@citelabel#1{b@#1}%
%
%
\def\@citedef#1#2{\expandafter\gdef\csname\@citelabel{#1}\endcsname{#2}}%
%
%
%
\def\@readbblfile{%
   \ifx\@itemnum\@undefined
      \@innernewcount\@itemnum
   \fi
   \begingroup
      \def\begin##1##2{%
         \setbox0 = \hbox{\biblabelcontents{##2}}%
         \biblabelwidth = \wd0
      }%
      \def\end##1{}
      %
      %
      \@itemnum = 0
      \def\bibitem{\@getoptionalarg\@bibitem}%
      \def\@bibitem{%
         \ifx\@optionalarg\empty
            \expandafter\@numberedbibitem
         \else
            \expandafter\@alphabibitem
         \fi
      }%
      \def\@alphabibitem##1{%
         \expandafter \xdef\csname\@citelabel{##1}\endcsname {\@optionalarg}%
         \ifx\biblabelprecontents\@undefined
            \let\biblabelprecontents = \relax
         \fi
         \ifx\biblabelpostcontents\@undefined
            \let\biblabelpostcontents = \hss
         \fi
         \@finishbibitem{##1}%
      }%
      \def\@numberedbibitem##1{%
         \advance\@itemnum by 1
         \expandafter \xdef\csname\@citelabel{##1}\endcsname{\number\@itemnum}%
         \ifx\biblabelprecontents\@undefined
            \let\biblabelprecontents = \hss
         \fi
         \ifx\biblabelpostcontents\@undefined
            \let\biblabelpostcontents = \relax
         \fi
         \@finishbibitem{##1}%
      }%
      \def\@finishbibitem##1{%
         \biblabelprint{\csname\@citelabel{##1}\endcsname}%
         \@writeaux{\string\@citedef{##1}{\csname\@citelabel{##1}\endcsname}}%
         \ignorespaces
      }%
      %
      %
      \let\em = \bblem
      \let\newblock = \bblnewblock
      \let\sc = \bblsc
      \frenchspacing
      \clubpenalty = 4000 \widowpenalty = 4000
      \tolerance = 10000 \hfuzz = .5pt
      \everypar = {\hangindent = \biblabelwidth
                      \advance\hangindent by \biblabelextraspace}%
      \bblrm
      \parskip = 1.5ex plus .5ex minus .5ex
      \biblabelextraspace = .5em
      \bblhook
      \input \bblfilebasename.bbl
   \endgroup
}%
%
%
\@innernewdimen\biblabelwidth
\@innernewdimen\biblabelextraspace
%
%
%
\def\biblabelprint#1{%
   \noindent
   \hbox to \biblabelwidth{%
      \biblabelprecontents
      \biblabelcontents{#1}%
      \biblabelpostcontents
   }%
   \kern\biblabelextraspace
}%
%
%
%
\def\biblabelcontents#1{{\bblrm [#1]}}%
%
%
\def\bblrm{\rm}%
%
%
\def\bblem{\it}%
%
%
\def\bblsc{\ifx\@scfont\@undefined
              \font\@scfont = cmcsc10
           \fi
           \@scfont
}%
%
%
\def\bblnewblock{\hskip .11em plus .33em minus .07em }%
%
%
\let\bblhook = \empty
%
%
%
\def\printcitestart{[}
\def\printcitefinish{]}
\def\printbetweencitations{, }
\def\printcitenote#1{, #1}
%
%
%
\let\citation = \@gobble
%
%
%
\@innernewcount\@numparams
%
%
\def\newcommand#1{%
   \def\@commandname{#1}%
   \@getoptionalarg\@continuenewcommand
}%
%
%
\def\@continuenewcommand{%
   \@numparams = \ifx\@optionalarg\empty 0\else\@optionalarg \fi \relax
   \@newcommand
}%
%
%
\def\@newcommand#1{%
   \def\@startdef{\expandafter\edef\@commandname}%
   \ifnum\@numparams=0
      \let\@paramdef = \empty
   \else
      \ifnum\@numparams>9
         \errmessage{\the\@numparams\space is too many parameters}%
      \else
         \ifnum\@numparams<0
            \errmessage{\the\@numparams\space is too few parameters}%
         \else
            \edef\@paramdef{%
               \ifcase\@numparams
                  \empty  No arguments.
               \or ####1%
               \or ####1####2%
               \or ####1####2####3%
               \or ####1####2####3####4%
               \or ####1####2####3####4####5%
               \or ####1####2####3####4####5####6%
               \or ####1####2####3####4####5####6####7%
               \or ####1####2####3####4####5####6####7####8%
               \or ####1####2####3####4####5####6####7####8####9%
               \fi
            }%
         \fi
      \fi
   \fi
   \expandafter\@startdef\@paramdef{#1}%
}%
%
%
%
%
\def\@readauxfile{%
   \if@auxfiledone \else 
      \global\@auxfiledonetrue
      \@testfileexistence{aux}%
      \if@fileexists
         \begingroup
            \endlinechar = -1
            \catcode`@ = 11
            \input \jobname.aux
         \endgroup
      \else
         \message{\@undefinedmessage}%
         \global\@citewarningfalse
      \fi
      \immediate\openout\@auxfile = \jobname.aux
   \fi
}%
%
%
\newif\if@auxfiledone
\ifx\noauxfile\@undefined \else \@auxfiledonetrue\fi
%
%
%
%
\@innernewwrite\@auxfile
\def\@writeaux#1{\ifx\noauxfile\@undefined \write\@auxfile{#1}\fi}%
%
%
%
\ifx\@undefinedmessage\@undefined
   \def\@undefinedmessage{No .aux file; I won't give you warnings about
                          undefined citations.}%
\fi
%
%
\@innernewif\if@citewarning
\ifx\noauxfile\@undefined \@citewarningtrue\fi
%
%
%
\catcode`@ = \@oldatcatcode


\def\widestnumber#1#2{}

\def\rm{\fam0 \tenrm}

\def\fakesubhead#1\endsubhead{\bigskip\noindent{\bf#1}\par}


%
%
%

%

\font\textrsfs=rsfs10
\font\scriptrsfs=rsfs7
\font\scriptscriptrsfs=rsfs5

\newfam\rsfsfam
\textfont\rsfsfam=\textrsfs
\scriptfont\rsfsfam=\scriptrsfs
\scriptscriptfont\rsfsfam=\scriptscriptrsfs

\edef\oldcatcodeofat{\the\catcode`\@}
\catcode`\@11

\def\Cal@@#1{\noaccents@ \fam \rsfsfam #1}

\catcode`\@\oldcatcodeofat

\bn

\newpage

\head {Anotated Content} \endhead  \resetall
\bn
\S0 $\quad$ Introduction
\mn
\S1 $\quad$ Countable Groups
\mr
\item "{{}}"  [We present a result on $\omega$ analytic equivalence relations
on ${\Cal P}(\omega)$ and apply it to $\aleph_0$-system of groups,
sharpening somewhat \cite{GrSh:302a}.]
\ermn
\S2 $\quad$ On $\lambda$-analytic equivalence relations
\mr
\item "{{}}"  [We generalize \S1 replacing $\aleph_0$ by $\lambda$ regular,
unfortunately this is only consistent.]
\ermn
\S3 $\quad$ On $\lambda$-systems of groups
\mr
\item "{{}}"  [This relates to \S2 as the application relates to the lemma
in \S1.]
\ermn
\S4 $\quad$ Back to the $p$-rank of Ext
\mr
\item "{{}}"  [We show that we can put the problem in the title to the
previous context, and show that in Easton model, \S2 and \S3 apply to every
regular $\lambda$.]
\ermn
\S5 $\quad$ Strong limit of countable cofinality
\mr
\item "{{}}"  [We continue \cite{GrSh:302a}.]
\endroster
\newpage

\head {\S0} \endhead  \resetall
\bn
A usual dichotomy 
is that in many cases, reasonably definable sets, satisfies
the CH, i.e. if they are uncountable they have cardinality continuum.  A 
strong dichotomy is when: if the cardinality is infinite it is continuum as in
\cite{Sh:273}.  We are interested in such phenomena when $\lambda =
\aleph_0$ is replaced by $\lambda$ regular uncountable and also by
$\lambda = \beth_\omega$ or more generally by 
strong limit of cofinality $\aleph_0$.
\bn
\ub{Question}:  Dom the parallel of \scite{1.2} holds for e.g.
$\beth_\omega$?  portion?
\newpage

\head {\S1 Countable groups} \endhead  \resetall
\bigskip

Here we give a complete proof a strengthening of the theorem of
\cite{GrSh:302a}, for the case $\lambda = \aleph_0$ using a variant of
\cite{Sh:273}.
\proclaim{\stag{1.1} Theorem}  Suppose
\mr
\item "{$(A)$}"  $\lambda$ is $\aleph_0$.  Let
$\langle G_m,\pi_{m,n}:m \le n < \omega \rangle$ be an 
inverse system whose inverse limit is $G_\omega$ with $\pi_{n,\omega}$ 
such that $|G_n| < \lambda$.  (So $\pi_{m,n}$ is
a homomorphism from $G_n$ to $G_m,\alpha \le \beta \le \gamma \le \omega
\Rightarrow \pi_{\alpha,\beta} \circ \pi_{\beta,\gamma} = 
\pi_{\alpha,\beta}$ and $\Pi_{\alpha,\alpha}$ is the identity).
\sn
\item "{$(B)$}"  Let $\bold I$ be an index set.  For every $t \in \bold I$,
let $\langle H^t_m,\pi^t_{m,n}:m \le n < \omega \rangle$ be an inverse system
of groups and $H^t_\omega$ with $\pi^t_{n,\omega}$ be the corresponding
inverse limit and $H^t_m$ of cardinality $< \lambda$.
\sn
\item "{$(C)$}"  Let for every $t \in \bold I,\sigma^t_n:H^t_n \rightarrow
G_n$ be a homomorphism such that all diagrams commute (i.e.
$\Pi_{m,n} \circ \sigma^t_n = \sigma^t_m \circ \pi^t_{m,n}$ for $m \le n <
\omega$), and let $\sigma^t_\omega$ be the induced homomorphism from $H^t_
\omega$ into $G_\omega$.
\sn
\item "{$(D)$}"  $\bold I$ is countable \footnote{this is stronger, earlier
$\bold I$ was finite}
\sn
\item "{$(E)$}"   For every $\mu < \lambda$ there is a sequence
$\langle f_i \in G_\omega:i < \mu \rangle$ such that for \nl
$i \ne j \and t \in \bold I \Rightarrow f_i f^{-1}_j \notin 
\text{ Rang}(\sigma^t_\omega)$. 
\ermn
\ub{Then} there is $\langle f_i \in G_\omega:i < 2^\lambda \rangle$ such
that \nl
$i \ne j \and t \in \bold I \Rightarrow f_i f^{-1}_j \notin 
\text{ Rang}(\sigma^t_\omega)$.
\endproclaim
\bn
This follows immediately from \scite{1.2} below.
\proclaim{\stag{1.2} Lemma}  Assume for every $n < \omega,{\Cal E}_n$ is an
analytic equivalence relation on ${\Cal P}(\omega) = \{A:A \subseteq 
\omega^+\}$ which satisfies
\mr
\item "{$(*)$}"  if $A,B \subset \bold Z^+,n \notin B,A = B \cup \{n\}$,
then $A,B$ are not ${\Cal E}_n$-equivalent.
\ermn
\ub{Then} there is a perfect subset of
${\Cal P}(\omega)$ of pairwise ${\Cal E}_n$-nonequivalent $A \subseteq
\omega$, simultaneously for all $n$.
\endproclaim
\bigskip

\remark{Remark}  The proof uses some knowledge of set theory and is close
to \cite[Lemma 1.3]{Sh:273}.
\endremark
\bigskip

\demo{Proof}  Let $r_m \in {}^\omega 2$ be the real parameter involved in
the definition $\varphi_m(x,y,r_m)$ of ${\Cal E}_m$.  Let $\bar \varphi =
\langle \varphi_m:m < \omega \rangle,\bar r = \langle r_m:m < \omega \rangle,
\bar{\Cal E} = \langle \bar{\Cal E}_m:m < \omega \rangle$.  Let $N$ be a
countable elementary submodel of $({\Cal H}((2^{\aleph_0})^+),\in)$
to which $\bar \varphi,\bar r,\bar{\Cal E}$ belong.  Now
\mr
\item "{$(**)$}"  if $\langle A_1,A_2 \rangle$ be a pair of subsets of
$\omega^+$ which is Cohen generic over $N$ [this means that it belongs to
no first category subset of ${\Cal P}(\omega) \times {\Cal P}(\omega)$ 
which belongs to $N$]
then
{\roster
\itemitem{ $(\alpha)$ }  $A_1,A_2$ are ${\Cal E}_m$-equivalent in
$N[A_1,A_2]$ if they are ${\Cal E}_m$-equivalent
\sn
\itemitem{ $(\beta)$ }   $A_1,A_2$ are non-${\Cal E}_m$-equivalent in
$N[A_1,A_2]$.
\endroster}
\endroster 
\enddemo
\bigskip

\demo{Proof of $(**)$}
\mr
\item "{$(\alpha)$}"  by the absoluteness criterions (Levy She\"onfied)
\sn
\item "{$(\beta)$}"  if not, then some finite information forces this,
hence for some $n$
{\roster
\itemitem{ $(*)$ }  if $\langle A'_1,A'_2 \rangle$ is Cohen generic over
$N$ and $A'_1 \cap \{0,1,\dotsc,n\} = A_1 \cap \{0,1,\dotsc,n\}$ and
$A'_2 \cap \{0,1,\dotsc,n\} = A_2 \cap \{1,\dotsc,n\}$ then $A'_1,A'_a2$
are ${\Cal E}_m$-equivalent in $N[A'_1,A'_2]$.
\endroster}
\ermn
Let $A''_1$ be $A_1 \cup\{n+1\}$ if $(n+1) \notin A_1$ and $A_1 \backslash
\{n+1\}$ if $(n+1) \in A_1$.

Trivially also $\langle A''_1,A_2 \rangle$ is Cohen generic over $N$, hence
by $(*)$ above $A''_1,A_2$ are ${\Cal E}_m$-equivalent in $N[A''_1,A_2]$.  By
$(**)(1)$ we know that really $A''_1,A_2$ are ${\Cal E}_m$-equivalent. 
As equivalence is a transitive relation clearly $A_1,A''_1$ are
${\Cal E}_m$-equivalent.   But this contradicts the hypothesis $(*)$.

We can easily find a perfect (nonempty) subset $P$ of $\{A:A \subseteq
\omega\}$ such that for any distinct $A,B \in P,(A,B)$ is Cohen generic over
$N$.  So for each $m$ for $ \ne B \in P,N[A,B] \models ``A,B$ are not
${\Cal E}_m$-equivalent" and by $(**)(\alpha)$ \, $A,B$ are not
${\Cal E}_m$-equivalent.  This finishes the proof. 
\hfill$\square_{\scite{1.2}}$
\enddemo
\bn
\centerline {$* \qquad * \qquad *$}
\bn
The same proof gives
\proclaim{\stag{1.3} Claim}  If ${\Cal E}$ is a two-place relation on
${\Cal P}(\omega)$ such that $A \subseteq \omega,B \subseteq \omega,A = B
\cup \{n\},n \notin B,C \subseteq \omega \Rightarrow \neg A {\Cal E} C \vee
\neg B \in {\Cal C}$, \ub{then} there is a perfect subset $\bold P$ of
${\Cal P}(\omega)$ such that $A \ne B \in \bold P \Rightarrow \neg A {\Cal E}
B$.
\endproclaim
\newpage

\head {\S2 On $\lambda$-analytic equivalence relations} \endhead  \resetall
\bigskip

\demo{\stag{2.0} Hypothesis}  $\lambda = \text{ cf}(\lambda)$ is fixed.
\enddemo
\bigskip

\definition{\stag{2.1} Definition}  1) A sequence $A$ of relations $\bar R =
\langle R_\varepsilon:\varepsilon < \varepsilon(*) \rangle$ on ${}^\lambda 2$
(equivalently ${\Cal P}(\lambda)$) i.e. a sequence of definitions of such
relations in $({\Cal H}(\lambda^+),\in)$ and with parameters in ${\Cal H}
(\lambda^+)$ is $\lambda$-w.c.a. (weakly Cohen absolute) if: for any $A 
\subseteq \lambda$
\mr
\item "{$(a)$}"  there are $N,r$ such that:
{\roster
\itemitem{ $(\alpha)$ }  $N$ is a transitive model
\sn
\itemitem{ $(\beta)$ }  $N^{< \lambda} \subseteq N,\lambda +1 \subseteq N$,
the sequence of the definitions of $\bar R$ belongs to $N$
\sn
\itemitem{ $(\gamma)$ }  $A \in N$
\sn
\itemitem{ $(\delta)$ }  $r \in {}^\lambda 2$ is Cohen over $N$ there is
generic for $({}^{\lambda>}2,\triangleleft)$ over $N$)
\sn
\itemitem{ $(\varepsilon)$ }  $R$ and $\neg R$ are absolute from $N[r]$ to
$V$.
\endroster}
\ermn
2) We say $R$ is $(\lambda,\mu)$-w.c.a. if for $A \subseteq \lambda$ we can
find $N,r_\alpha$ (for $\alpha \in 2^\lambda$) such that satisfying clauses
$(\alpha),(\beta),(\gamma)$ from above and
\mr
\item "{$(\delta)'$}"  for $\eta \ne \nu \in {}^\lambda 2,(\eta,\nu)$ is a
pair of Cohens over $N$
\sn
\item "{$(\varepsilon)'$}"  $R,\neg R$ are absolute from $N[\eta,\nu]$ to
$V$.
\ermn
3) We say $\lambda$ is $(\lambda,\mu)$-w.c.a. \ub{if} every $\lambda$-analytic
relation on ${}^\lambda 2$ is $(\lambda,\mu)$-w.c.a. \nl
That is, it has the form $(\exists Y_1,\dotsc,Y_m \subseteq \lambda \times
\lambda)\varphi(Y_1,\dotsc,Y_m;X_1,\dotsc,X_n)$.
\enddefinition
\bigskip

\proclaim{\stag{2.2} Claim}  Assume
\mr
\item "{$(a)$}"  for $\varepsilon < \varepsilon(*) \le \lambda,
{\Cal E}_\varepsilon$ is a $(\lambda,\mu)$-w.c.a. equivalence relation on
${\Cal P}(\lambda)$, more exactly a definition of one and
\sn
\item "{$(b)$}"  if $A,B \subseteq \lambda$ and $\alpha \in A \backslash B,
A = B \cup \{\alpha\}$, \ub{then} $A,B$ are not ${\Cal E}$-equivalent.
\ermn
\ub{Then} there is a set ${\Cal P} \subseteq {\Cal P}(\lambda)$ of
$\mu$-pairwise non-${\Cal E}_\varepsilon$-equivalent members of ${\Cal P}
(\lambda)$ for all $\varepsilon < \varepsilon(*)$ simultaneously.
\endproclaim
\bigskip

\remark{\stag{2.3} Remark}  If in \scite{2.1} we ask that $\{r_\eta:\eta \in
{}^\lambda 2\}$ perfect, then so is ${\Cal P}$.
\endremark
\bigskip

\definition{\stag{2.4} Definition}  1) ${\Cal P} \subseteq {\Cal P}(\lambda)$
is perfect if there is a $\lambda$-perfect $T \subseteq {}^{\lambda >} 2$
(see below) such that ${\Cal P} = \{\{\alpha < \lambda:\eta(\alpha) =1\}:
\eta \in \text{ lim}_\lambda(T)\}$. \nl
2) $T$ is a $\lambda$-perfect is:
\mr
\item "{$(a)$}"  $T \subseteq {}^{\lambda >}2$ is non-empty
\sn
\item "{$(b)$}"  $\eta \in T \and \alpha < \ell g(\eta) \Rightarrow \eta
\restriction \alpha \in T$
\sn
\item "{$(c)$}"  if 
$\delta < \lambda$ is a limit ordinal, $\eta \in {}^\delta 2$
and $(\forall \alpha < \delta)(\eta \restriction \alpha \in T)$, \ub{then}
$\eta \in T$
\sn
\item "{$(d)$}"  if $\eta \in T,\ell g(\eta) < \alpha < \lambda$ then there
is $\nu,\eta \triangleleft \nu \in T \cap {}^\alpha 2$
\sn
\item "{$(e)$}"  if $\eta \in T$ then there are $\triangleleft$-incomparable
$\nu_1,\nu_2 \in T$ such that \nl
$\eta \triangleleft \nu_1 \and \eta \triangleleft \nu_1$.
\ermn
3) Lim$_\delta(T) = \{\eta:\ell g(\eta) = \delta$ and $(\forall \alpha <
\delta)(\eta \restriction \alpha \in T)\}$.
\enddefinition
\bigskip

\demo{Proof of \scite{2.2}}

Let $T^* = {}^{\lambda >}2$. \nl
Let $N$ and $r_\alpha \in {}^\lambda 2$ for $\alpha < \mu$ be as in
Definition \scite{2.1}.  We identify $r_\alpha$ with $\{\gamma < \lambda:
r_\alpha(\gamma) = 1\}$.

It is enough to prove that, assuming $\alpha \ne \beta < \mu$
\mr
\item "{$(*)_1$}"  $\neg r_\alpha {\Cal E}_\varepsilon r_\beta$.
\ermn
By clause $(\varepsilon)$ of Definition \scite{2.1}(2) it is enough to prove
\mr
\item "{$(*)_2$}"  $N[r_\alpha,r_\beta] \ne \neg r_\alpha 
{\Cal E}_\varepsilon \nu_\beta$.
\ermn
Assume this fails, so $N[r_\alpha,r_\beta] \ne r_\alpha 
{\Cal E}_\varepsilon r_\beta$ then for some $i< \lambda$

$$
(r_\alpha \restriction i,r_\beta \restriction i) 
\Vdash_{({}^{\lambda >}2) \times
({}^{\lambda >}2)} ``{\underset\tilde {}\to r_1} {\Cal E}_\varepsilon
{\underset\tilde {}\to r_2}".
$$
\mn
Define $r \in {}^\lambda 2$ by

$$
r(j) = \cases r_\beta(j) &\text{ \ub{if} } j \ne i \\
1-r_\beta(j) &\text{ \ub{if} } j=i \endcases
$$
\mn
So also $(r_\alpha,r)$ is a generic pair for ${}^{\lambda >}2 \times
{}^{\lambda >}2$ over $N$ and $(r_\alpha \restriction i,r \restriction i) 
=$ \nl
$(r_\alpha \restriction i,r_\beta \restriction i)$ 
hence by the forcing theorem

$$
N[r_\alpha,r] \models {\underset\tilde {}\to r_\alpha} 
{\Cal E}_\varepsilon r.
$$
\mn
By $r_\alpha,r_\beta,r \in N[r_\alpha,r_\beta] = N[r_\alpha,r]$, so
$N[r_\alpha,r_\beta] \models r_\beta {\Cal E}_\varepsilon r$ hence $V \models 
r_\beta {\Cal E}_\varepsilon r$, a contradiction to assumption (b).
\hfill$\square_{\scite{2.2}}$
\enddemo
\bigskip

\definition{\stag{2.5} Definition}  We call $Q$ a pseudo $\lambda$-Cohen
forcing if:
\mr
\item "{$(a)$}"  $Q = \{p:p \text{ a partial functiom from } \lambda
\text{ to } \{0,1\}\}$
\sn
\item "{$(b)$}"  $p \le_Q q \Rightarrow p \subseteq q$
\sn
\item "{$(c)$}"  ${\Cal I}_i = \{p:i \in \text{ Dom}(p)\}$ is dense for $i <
\lambda$
\sn
\item "{$(d)$}"  define $F_i:{\Cal I}_i \rightarrow {\Cal I}_i$ by: 
Dom$(F_i(p)) = \text{ Dom}(F_i(p))$

$$
(F_i(p))(j) = \cases p(j) &\text{ \ub{if} } j=i \\
1-p(j) &\text{ \ub{if} }j \ne i \endcases
$$
\mn
then $F_i$ is an automorphism of $({\Cal I}_i,<^Q \restriction {\Cal I}_i)$.
\endroster
\enddefinition
\bigskip

\proclaim{\stag{2.6} Claim}  In \scite{2.1}, \scite{2.4} we can replace
$({}^{\lambda >}2,\triangleleft)$ by $Q$.
\endproclaim
\bn
\ub{\stag{2.7} Observation}:  So in $V \models G.C.H., P$ is Easton forcing,
\ub{then} in $V^P$ for every $\lambda$ regular for $Q = (({}^{\lambda>}2)^V,
\triangleleft)$
we have: $Q$ is pseudo $\lambda$-Cohen and in $V^P$ we have $\lambda$ is
$(\lambda,2^\lambda)$-w.c.a.
\bn
\ub{\stag{2.8} Discussion}:  But in fact $\lambda$ being $(\lambda,2^\lambda)$-
w.c.a. is a weak condition.
\bn
We can generalize further using the following definition
\definition{\stag{2.10} Definition}  1) For $r_0,r_1 \in {}^\lambda 2$ we say
$(r_0,r_1)$ or $r_0,r_1$ is an $\bar R$-pseudo Cohen pair over $N$ if 
($\bar R$ is a definition (in $({\Cal H}(\lambda^+),\in))$ 
of a relation on ${\Cal P}(\lambda)$ (or ${}^\lambda 2$) and)
for some forcing notion $Q$ and $Q$-names ${\underset\tilde {}\to r_0},
{\underset\tilde {}\to r_1}$ and $G \subseteq Q \, (G \in V)$ generic over
$N$ we have:
\mr
\item "{$(a)$}"  ${\underset\tilde {}\to r_0}[G] = r_0,
{\underset\tilde {}\to r_1}[G] = r_1$
\sn
\item "{$(b)$}"  for every $p \in G,i < \lambda$ large enough and $\ell(*)
< 2$ there is $G' \subseteq Q$ generic over $N$ such that: $p \in G$ and
$({\underset\tilde {}\to r_\ell}[G'])(j) = ({\underset\tilde {}\to r_\ell}
[G])(j) \Leftrightarrow (j,\ell) \ne (i,\ell(*))$
\sn
\item "{$(c)$}"  $R$ is absolute from $N[G]$ and from $N[G']$ to $V$.
\ermn
2) We say $\lambda$ is $\mu$-p.c.a for $\bar R$ \ub{if} for every $x \in
{\Cal H}(\lambda^+)$ there are $N, \langle \eta_i:i < \mu \rangle$ such that:
\mr
\item "{$(a)$}"  $N$ is a transitive model of $ZFC^-$
\sn
\item "{$(b)$}"  for $i \ne j < \mu, (r_i,r_j)$ is an $\bar R$-pseudo Cohen
pair over $N$.
\ermn
3) We omit $R$ if this holds for any $\lambda$-sequence of $\sum^1_1$
formula in ${\Cal H}(\lambda^+)$.
\enddefinition
\bn
Clearly
\proclaim{\stag{2.11} Claim}  1) If $\lambda$ is $\mu$-p.c.a for ${\Cal E},
{\Cal E}$ an equivalence relation on ${\Cal P}(\lambda)$ and $A \subseteq B
\subseteq \lambda \and |B \backslash A| = 1 \Rightarrow \neg A {\Cal E} B$,
\ub{then} ${\Cal E}$ has $\ge \mu$ equivalence classes. \nl
2) Similarly if ${\Cal E} = \dsize \bigvee_{\varepsilon < \varepsilon(*)}
{\Cal E}_\varepsilon,\varepsilon(*) \le \lambda$ and each 
${\Cal E}_\varepsilon$ as above.
\endproclaim
\newpage

\head {\S3 On $\lambda$-systems of groups} \endhead  \resetall
\bigskip

\demo{\stag{3.0} Hypothesis}  $\lambda = \text{ cf}(\lambda)$.

We may wonder does \scite{3.2} have any cases it covers?
\enddemo
\bigskip

\definition{\stag{3.1} Definition}  1) We say ${\Cal Y} = (\bar A, \bar K,
\bar G,\bar D,g^*)$ is a $\lambda$-system if
\mr
\item "{$(A)$}"  $\bar A = \langle A_i:i \le \lambda \rangle$ is an
increasing sequence of sets, $A = A_\lambda = \dbcu_i A_i$
\sn
\item "{$(B)$}"  $\bar K = \langle K_t:t \in A \rangle$ is a sequence of
finite groups
\sn
\item "{$(C)$}" $\bar G = \langle G_i:i \le \lambda \rangle$ is a sequence
of groups, $G_i \subseteq \dsize \prod_{t \in A_i} K_t$, each $G_i$ is
closed and $i < j \le \lambda \Rightarrow G_i = \{g \restriction A_i:g \in
G_j\}$ and \nl
$G_\lambda = \{g \in \dsize \prod_{t \in A_\lambda} K_t:(\forall
i < \lambda)(g \restriction A_i \in G_i)$
\sn
\item "{$(D)$}"  $\bar D = \langle D_\delta:\delta \le \lambda \text{ (a limit
ordinal) } \rangle,D_\delta$ an ultrafilter on $\delta$ such that \nl
$\alpha < \delta \Rightarrow [\alpha,\delta) \in D_\delta$
\sn
\item "{$(E)$}"  $\bar g^* = \langle g_i:i < \lambda  \rangle,g_i \in
G_\lambda$ and $g_i \restriction A_i = e_{G_i} = \langle e_{K_t}:t \in A_i
\rangle$.
\ermn
2) Let ${\Cal Y}^-$ be the same omitting $\bar D$.
\enddefinition
\bigskip

\definition{\stag{3.2} Definition}  For a $\lambda$-system ${\Cal Y}$ and
$j \le \lambda + 1$ we say $\bar f_i \in \text{ cont}(j,{\Cal Y})$ if:
\mr
\item "{$(a)$}"  $\bar f = \langle f_i:i < j \rangle$
\sn
\item "{$(b)$}"  $f_i \in G_\lambda$
\sn
\item "{$(c)$}"  if $\delta < j$ is a limit ordinal then $f_\delta =
\text{ Lim}_{D_\delta}(\bar f \restriction \delta)$ which means:
$$
\text{ for every } t \in A,f_\delta(t) = \text{ Lim}_{D_\delta} \langle f_i
(t):i < \delta \rangle
$$
\sn
which means
$$
\{i < \delta:f_\delta(t) = f_i(t)\} \in D_\delta.
$$
\endroster
\enddefinition
\bn
\ub{\stag{3.3} Fact}:  1)  If $\bar f \in \text{ cont}(j,{\Cal Y}),i < j$ 
\ub{then} $\bar f \restriction i \in \text{ cont}(i,{\Cal Y})$. \nl
2) If $\bar f \in \text{ cont}(j,{\Cal Y})$ and $j < \lambda$ is non-limit,
and $f_j \in G_\lambda$ \ub{then}

$$
\bar f \char 94 \langle f_j \rangle \in \text{ cont}(j+1,{\Cal Y}).
$$
\mn
3) If $\bar f \in \text{ cont}(j,{\Cal Y})$ and $j$ is a limit ordinal
$\le \lambda$, \ub{then} for some unique $f_j \in G_\lambda$ we have
$\bar f \char 94 \langle f_j \rangle \in \text{ cont}(j+1,{\Cal Y})$. \nl
4) If $j \le \lambda +1,f \in G$ \ub{then} $\bar f = \langle f:i < j \rangle
\in \text{ cont}(j,{\Cal Y})$. \nl
5) If $\bar f,\bar g \in \text{ cont}(j,{\Cal Y})$, \ub{then}
$\langle f_ig_i:i < j \rangle$ and $\langle f^{-1}_i:i <j \rangle$ belongs
to $\text{cont}(j,{\Cal Y})$.
\bigskip

\demo{Proof}  Straight.
\enddemo
\bigskip

\definition{\stag{3.4} Definition}  1) If $\bar g \in {}^j(G_\lambda),
j \le \lambda$ we define $f_{\bar g}$ by induction on $j$ for all such 
$\bar g$ as follows:
\mn

\ub{$j=0$}:  $f_{\bar g} = e_G = \langle e_{K_t}:t \in A \rangle$
\mn

\ub{$j=i+1$}:  $f_{\bar g} = f_{\bar g \restriction i} g_i$
\mn

\ub{$j$ limit}:  $f_{\bar g} = \text{ Lim}_{D_\delta} \langle 
f_{\bar g \restriction i}:i < j \rangle$
\mn
2) We say $\bar g$ is trivial on $X$ if $i \in X \cap \ell g(\bar g)
\Rightarrow g_i = e_{G_\lambda}$. \nl
3) For $\eta \in {}^{\lambda >}2$ let $\bar g^\eta = \langle f^\eta_i:i <
\ell g(\eta) \rangle$, where
$$
g^\eta_i = \cases g^*_i &\text{ if } \eta(i)=1 \\
  e_{G_\lambda} &\text{ if } \eta(i) = 0 \endcases
$$
\enddefinition
\bigskip

\proclaim{\stag{3.5} Claim} 1) If $i \le j$ and $\bar g,\bar g',\bar g'' \in
{}^j(G_\lambda),\bar g' \restriction i = \bar g \restriction i,\bar g'$
is trivial on $[i,j)$, \nl
$\bar g'' \restriction [i,j) = \bar g \restriction
[i,j)$ and $\bar g'''$ is trivial on $i$, \ub{then}:

$$
f_{\bar g} = f_{\bar g'} f_{\bar g''} \text{ and } f_{\bar g'} =
f_{\bar g \restriction i}.
$$
\mn
2) For $\eta \in {}^\lambda 2,f_{(\bar g^\eta)} = \text{ Lim}\langle
f_{(\bar g,\eta \restriction i)}:i < \lambda \rangle$ (i.e. any ultrafilter
$D'_\lambda$ on $\lambda$ containing the co-bounded sets will do).
\endproclaim
\bigskip

\demo{Proof}  Straight.
\enddemo
\bigskip

\proclaim{\stag{3.6} Claim}  1) Let ${\Cal Y}^-$ be a $\lambda$-system, 
$H_\varepsilon$ a subgroup of $G_\lambda$, for $\varepsilon < \varepsilon(*) 
< \lambda$ and ${\Cal E}_\varepsilon$ the equivalence 
relation $\dsize \bigwedge_{t \in \bold I} [f'(f'')^{-1} \in H_t]$ 
and assume: $\lambda > i \ge \varepsilon \Rightarrow g^*_i 
\notin H_\varepsilon$
\mr
\item "{$(\alpha)$}"  the assumption $(A)$ of \scite{3.2} holds with
$f_A = f_{(\bar g^\eta)}$ when $A \subseteq \lambda,\eta \in {}^\lambda 2,
A = \{i:\eta(i) = 1\}$
\sn
\item "{$(\beta)$}"  if in addition $A_\varepsilon \subseteq \lambda^+,
K_t \in {\Cal H}(\lambda^+)$ and $\langle H_\varepsilon:\varepsilon <
\varepsilon(*) \rangle$ is $(\lambda,\mu)$-w.c.a., \ub{then} also assumption
$(B)$ of \scite{3.2} holds (hence its conclusion).
\endroster
\endproclaim
\bigskip

\proclaim{\stag{3.6A} Claim}  Assume
\mr
\item "{$(A)$}"  ${\Cal Y}^-$ a $\lambda$-system, $A_i \subseteq \lambda^+,
|A_i| \le \lambda,G_i \in {\Cal H}(\lambda^+)$
\sn
\item "{$(B)$}"  $\varepsilon(*) \le \lambda,\bar H = \langle H^\varepsilon_i:
i \le i(*) \rangle,\Pi^\varepsilon_{i,j}:H^\varepsilon_j \rightarrow
H^\varepsilon_i$ a homomorphism, for $i_0 \le i_1 \le i_2$ we have
$\Pi^\varepsilon_{i_0,i_1} \circ \Pi^\varepsilon_{i_1,i_2} = \Pi^\varepsilon
_{i_0,i_2},\sigma^\varepsilon_i:H^\varepsilon_t \rightarrow G_i,
\sigma^\varepsilon_i \Pi^\varepsilon_{i,j}(f) = (\sigma^\varepsilon_j(f))
\restriction A_i,H^\varepsilon_\lambda,\sigma^\varepsilon_\lambda$ is the
inverse limit (with $\Pi^\varepsilon_{i,\lambda}$) of $\langle
H^\varepsilon_i,\sigma^\varepsilon_{i,j}:i \le j < \lambda \rangle$ and 
$i < \lambda \Rightarrow H^\varepsilon_i \in {\Cal H}(\lambda^+)$
\sn
\item "{$(C)$}"  $H_\varepsilon = \text{ Rang}(\sigma^\varepsilon_\lambda)$.
\ermn
Then
\mr
\item "{$(\alpha)$}"  the assumptions of \scite{4.6} holds
\sn
\item "{$(\beta)$}"  if $\lambda$ is $(\lambda,\mu)$-w.c.a. \ub{then} also the
conclusion of \scite{4.6}, \scite{3.2} holds.
\endroster
\endproclaim
\bigskip

\demo{Proof}  Straight.
\enddemo
\bigskip

\centerline {$* \qquad * \qquad *$}
\bn
We can go one more step in concretization.
\proclaim{\stag{3.7} Claim}  1) Assume
\mr
\item "{$(a)$}"  $L$ is an abelian group of cardinality $\lambda$
\sn
\item "{$(b)$}"  $p$ a prime number
\sn
\item "{$(c)$}"  if $L' \subseteq L,|L'| < \lambda$, then Ext$_p(L,L',
\Bbb Z) \ne 0$
\sn
\item "{$(d)$}"  $\lambda$ is $\mu$-w.c.a. (in $V$).
\ermn
\ub{Then} $\mu \le r_p$(Ext$(G,\Bbb Z))$. \nl
2) If $(a),(b),(d)$ above, $\mu > \lambda,\lambda$ strongly inaccessible
then $r_p(\text{Ext}(G,\Bbb Z)) \notin [\lambda,\mu)$.
\endproclaim
\bigskip

\demo{Proof}  Without loss of generality $G$ is $\aleph_1$-free (so torsion
free).
\nl
Without loss of generality the set of elements of $G$ is
$\lambda$.  Let $A = A_\lambda = \lambda,L_\lambda = L$, for $j < \lambda,
A_j$ a proper initial segment of $\lambda$ such that $L_j = L \restriction
A_j$ is a pure subgroup of $L$, increasing continuously with $j$.  \nl
Let $K_t = \Bbb Z/p \Bbb Z,G_i = f \in \text{ HOM}(L_i,\Bbb Z/p \Bbb Z)$. \nl
Let $\varepsilon(*) = 1,H_i = \text{ HOM}(L_i,\Bbb Z)$ and 
$(\sigma^\varepsilon_i(f))(x) = f(x) + p \Bbb Z,K_\varepsilon = \text{ Rang}
(\sigma^\varepsilon_\lambda)$.  We know that $r_p(\text{Ext}(G,\Bbb Z))$ is
$(G_\lambda:K_0)$.  By assumption $(d)$ for each $i < \lambda$ we
can choose $g^*_i \in G_\lambda \backslash K_\varepsilon$ such that
$g^*_i \restriction L_i$ is zero. \hfill$\square_{\scite{3.7}}$
\enddemo
\newpage

\head {\S4 Back to the $p$-rank of Ext} \endhead  \resetall
\bigskip

\definition{\stag{4.1} Definition}  Let

$$
\align
\Xi_{\Bbb Z} = \bigl\{ \bar \lambda:&\, \bar \lambda = \langle \lambda_p:p <
\omega \text{ prime or zero}\rangle \text{ and for some} \\
  &\text{ abelian } (\aleph_1 \text{-free) group } L,\lambda_p = r_p
\text{ (Ext} (G,\Bbb Z)) \bigr\}.
\endalign
$$
\mn
Clearly $\Xi_{\Bbb Z}$ is closed under products.  Let $\bold P$ be the set
of primes.
\enddefinition
\bn
Remember that (see \cite[AP]{Sh:f}, 2.7, 2.7A, 2.13(1),(2)). \nl
\ub{\stag{4.2} Fact}:  1) In the Easton model if $G$ is $\aleph_1$-free not
free, $G' \subseteq G,|G'| < |G| \Rightarrow G/G'$ not free \ub{then}
$r_0(\text{Ext}(G,\Bbb Z)) = 2^{|G|}$.
\bn
\ub{\stag{4.3} Fact}:  1) Assume $\mu$ is strong limit $> \aleph_0$,
cf$(\mu) = \aleph_0,\lambda = \mu,2^\mu = \mu^+$ and some 
$Y \subseteq [{}^\omega \mu]
^\lambda$ is $\mu$-free, (equivalently $\mu^+$-free, see in proof). \nl
Let $\bold P_0,\bold P_1$ be a partition of the primes. \nl
\ub{Then} for some $\aleph_1$-free abelian group $L,|L| = \mu^+,2^\lambda
=r_0$ (Ext$(G,\Bbb Z)),p \in \bold P_1 \Rightarrow \lambda_p(\text{Ext}
(G,\Bbb Z)) = 2^\lambda,p \in \bold P_0 \Rightarrow \lambda_p(\text{Ext}
(G,\Bbb Z)) = 0$. 
\bigskip

\remark{Remark}  On other cardinals see \cite{MRSh:314}.
\endremark
\bigskip

\demo{Proof}  Let $Y = \{\eta_i:i < \lambda\}$.  Let $r:\mu^2 \rightarrow
\mu$ be a pairing function, so $pr(pr_1(\alpha),pr_2(\alpha)) = \alpha$.
Without loss of generality $\eta_i(n) = \eta_j(m) \Rightarrow n=m \and
\eta_i \restriction m = \eta_j \restriction m$.  Let $L_0$ be
${\underset {\alpha < \lambda} {}\to \bigoplus} \Bbb Z x_\alpha$.  Let
$\langle f_i:i < \lambda \rangle$ list the functions $f \in \text{ HOM}
(L,\Bbb Z/p \Bbb Z)$.  We choose by induction on $i < \lambda,(g_i,\nu_i,
\ell_i)$ such that:
\mr
\item "{$(\alpha)$}"  $g_i \in \text{ HOM}(L,\Bbb Z)$
\sn
\item "{$(\beta)$}"  $(\forall x \in L)[g_i(x)/p \Bbb Z = f_i(x)]$
\sn
\item "{$(\gamma)$}"  $\rho_i,\nu_i \in {}^\omega \mu,\eta_i(n) = pr_1
(\nu_i(n)) = pr_1(\rho_i(n))$
\sn
\item "{$(\delta)$}"  $(\forall j \le i)(\exists n < \omega)(\forall n)
[n \le m < \omega \rightarrow g_j(x_{\nu_i(m)}) = g_j(x_{\rho_i(m)})$
\sn
\item {$(\varepsilon)$}"  $(\forall j < i)(\exists n < \omega)(\forall m)
(n \le m < \omega \rightarrow g_i(x_{\nu_j(m)}) = g_i(x_{\nu_j(m)})$
\sn
\item "{$(\zeta)$}"  $\nu_i(n) \ne \rho_i(m)$.
\ermn
Arriving to $i$ let $h_i:i \rightarrow \omega$ be such that $\langle \{
\eta_j \restriction \ell:\ell \in [h_i(j),\omega)\}:j < i\rangle$ are 
pairwise disjoint (possible as $Y$ is $\mu^+$-free).  Now choose $g_i$ such
that clause $(\varepsilon)$ holds with $n = h_i(j)$, as the choice of
$h,\sigma$ splits the problem.  Let $i = \dbcu_{n < \omega} A^i_n,A^i_n
\subseteq A^i_{n+1},|A^i_n| < \mu$.  Now choose by induction on $n,\alpha
= \rho_i(n),\beta = \nu_i(n)$ as distinct ordinals $\in \{ \alpha \in \mu:
\alpha \notin \{\nu_i(m),\rho_i(m:m < m\}$ and $pr_1(\alpha) = \eta_i(n)\}$
such that $\langle g_j(x_{\nu_\alpha}):j \in A^i_n \rangle = \langle g_j
(x_\beta):j \in A^i_n \rangle$. \nl
Let $L$ be generated by $L_0 \cup \{y_{i,m}:i < \lambda,m < \omega\}$ freely
except that (the equations of $L_0$ and) $(\dsize \prod_{p \in \bold P_0
\cap n} p) y_{i,n+1} = y_{i,n} + x_{\nu_i(n)} - x_{\rho_i(n)}$.
\enddemo
\bn
\ub{Discussion}:  Suppose $V \models 2^{\aleph_\alpha} = \aleph_{\alpha +2}$.
If we use the Easton forcing, $2^{\aleph_n} = \aleph_{n+2}$ we get:
$|G| < \aleph_\omega,G$ not free, $\lambda_p(\text{Ext}(G,\Bbb Z))$
uncountable two.  In $\aleph_{\omega+1},\aleph_{\omega +2}$ have to look at
\cite{GiSh:597}.  So CON$[(\forall L)(\text{Ext}(L,\Bbb Z) \ne 0 \rightarrow
(\text{Ent}_p(L,\Bbb Z) \ne 0)]$ is hard.
\bn
\ub{\stag{4.4} Question}:  Do we have compactness for singular for
Ext$_p(G,\Bbb Z) = 0$?
\bigskip

\proclaim{\stag{4.5} Claim}  Let $V \models \kappa$ is supercompact,
$\lambda = \text{ cf}(\lambda) > \kappa, GCH$ and $\lambda =
\text{ cf}(\lambda) > \lambda_0 \Rightarrow \diamondsuit^*_\lambda$. \nl
\ub{Then} in some generic extension
\mr
\item "{$(*)$}"  if $L$ is $\aleph_1$-free abelian group, not free 
$L' \subseteq L \and |L'| < L \Rightarrow L/L'$ not free then
{\roster
\itemitem{ $(\alpha)$ }  $\lambda_0(\text{Ext}(L,\Bbb Z)) = 2^{|L|}$
\sn
\itemitem{ $(\beta)$ }   $\lambda_p(\text{Ext}(L,\Bbb Z)) = 2^{|L|}$.
\endroster}
\endroster
\endproclaim
\bigskip

\demo{Proof}  Let $P$ be adding $\lambda$ Cohen reals.
\enddemo
\bn
\ub{\stag{4.10} Question}:  If $\bar \lambda \in \Xi_{\Bbb Z}$ can we derive
$\bar \lambda' \in \Xi_{\Bbb Z}$ by increasing some $\lambda_p$'s?
\bn
\ub{\stag{4.7} Question}:  Suppose we choose $\lambda > 2^{\aleph_0}$
strongly inaccessible, $p^*$ prime and add by iterated forcing a $G$
satisfying $r_0 (\text{Ext}(G,\Bbb Z)) = 2^\lambda,r_{p^*}(\text{Ext}
(G,\Bbb Z)) =0$.  Does
not $\bar \lambda \in \Xi_{\Bbb Z} \and p \ne p^* \Rightarrow \lambda_p =
\lambda_0$?
\bn
\ub{\stag{4.8} Fact}:  If $\bar \lambda^i = 
\langle \lambda^i_p:p \in \bold P \cup
\{0\} \rangle \in \Xi_{\Bbb Z}$ for $i < \alpha$ and $\lambda_p =
\dsize \prod_{i < \alpha} \lambda^i_p$, then \nl
$\langle \lambda_p:p \in \bold P \cup \{0\} \rangle \in \Xi_{\Bbb Z}$.
\bigskip

\demo{Proof}  As if $G = {\underset {i < \alpha} {}\to \bigoplus}$ then
Ext$(G,\Bbb Z) = \dsize \prod_{i < \lambda} \text{ Ext}(G,\Bbb Z)$ hence
$\lambda_p(\text{Ext}(G,\Bbb Z)) = \dsize \prod_{i < \alpha} \lambda_p
(\text{Ext}(G_i,\Bbb Z))$.
\enddemo
\bn
\ub{\stag{4.9} Concluding Remark}:  In 
\cite{EkSh:505} the statement ``there is a $W$-abelian group" is 
characterized.

We can similarly characterize ``there is a separable group".
We have the same characterization for ``there
is a non-free abelian group" such that for some $p$, \nl
$r_p(\text{Ext}(G,\Bbb Z))=0$.
\bn
\ub{Question}:  Can $\bold P^* = \{p:p$ prime and $\bar \lambda \in
\Xi_{\Bbb Z} \and \lambda_0 > 0 \Rightarrow \lambda_p > 0\}$ be
$\notin \{\emptyset,\bold P\}$?
\bigskip

\proclaim{\stag{4.6} Claim}  If 
$\lambda$ is strong inaccessible or $\lambda = \mu^+,
\mu$ strong limit singular of cofinality $\aleph_0,S \subseteq \{\delta <
\lambda:\text{cf}(\delta) = \aleph_0\}$ is stationary not reflecting and
$\diamondsuit^*_S$ and $\bold P_0$ a set of primes, \ub{then} there is 
a $\lambda$-free abelian group $G$ such that $r_0(\text{Ext}(G,\Bbb Z) = 
2^\lambda,r_p(\text{Ext}(G,\Bbb Z)) = 0$ and: 
$p \in \bold P_0 \Rightarrow r_p(\text{Ext}(G,\Bbb Z)) = 2^\lambda$ and
$p$ prime and $p \notin \bold P_0 \Rightarrow r_p(\text{Ext}(G,\Bbb Z)
=0$.
\endproclaim
\newpage

\head {\S5 Strong limit of countable cofinality} \endhead  \resetall
\bn
We continue \cite{GrSh:302a}.  For consistency of ``no examples" see
\cite{MkSh:418}.
\definition{\stag{5.1} Theorem}  1) We say ${\Cal A}$ is a $(\lambda,
\bold I)$-system if ${\Cal A} = (\lambda,\bold I,\bar G = \langle
G_\alpha:\alpha \le \omega \rangle,
\bar H^t = \langle H^t_\alpha:\alpha \le \omega \rangle,
\bar \pi = \langle \pi_{\alpha,\beta},
\pi^t_{\alpha,\beta}:\alpha \le \beta \le \omega,t \in \bold I \rangle,
\bar \sigma = \langle \sigma_t:t \in \bold I \rangle)$ satisfies 
(we may write $\lambda^{\Cal A},
\pi^{t,{\Cal A}}_{\alpha,\beta}$, etc.)
\mr
\item "{$(A)$}"  $\lambda$ is $\aleph_0$ or generally a 
cardinal of cofinality $\aleph_0$  
\sn
\item "{$(B)$}"  $\langle G_m,\pi_{m,n}:m \le n < \omega \rangle$ is an 
inverse system whose inverse limit is $G_\omega$ with $\pi_{n,\omega}$ 
such that $|G_n| \le \lambda$.  (So $\pi_{m,n}$ is
a homomorphism from $G_n$ to $G_m,\alpha \le \beta \le \gamma \le \omega
\Rightarrow \pi_{\alpha,\beta} \circ \pi_{\beta,\gamma} = \pi_{\alpha,\beta}$
and $\pi_{\alpha,\alpha}$ is the identity).
\sn
\item "{$(C)$}"  $\bold I$ is an index set of cardinality $\le \lambda$. 
For every $t \in \bold I$ we have \nl
$\langle H^t_m,\pi^t_{m,n}:m \le n < \omega \rangle$ is an inverse system
of groups and $H^t_\omega$ with $\pi^t_{n,\omega}$ being the corresponding
inverse limit $H^t_\omega$ with $\pi^t_{m,\omega}$ and $H^t_m$ has 
cardinality $\le \lambda$.
\sn
\item "{$(D)$}"  for every $t \in \bold I,\sigma^t_n:H^t_n \rightarrow
G_n$ is a homomorphism such that all diagrams commute (i.e.
$\pi_{m,n} \circ \sigma^t_n = \sigma^t_m \circ \pi^t_{m,n}$ for $m \le n <
\omega$), and let $\sigma^t_\omega$ be the induced homomorphism from $H^t_
\omega$ into $G_\omega$
\sn
\item "{$(E)$}"  $G_0 = \{e_{G_0}\},H^t_0 = \{e_{H^t_0}\}$ (just for 
simplicity).
\ermn
2) We say ${\Cal A}$ is strict if 
$|G_n| < \lambda,|H^t_n| < \lambda,|\bold I| < \lambda$.  Let 
${\Cal E}_t$ be the following equivalence relation on $G_\omega:
f {\Cal E}_t g$ iff $fg^{-1} \in \text{ Rang}(\sigma^t_\omega)$. \nl
3) Let nu$({\Cal A}) = \sup\{\mu:\text{for each } n < \omega,
\text{ there is a sequence } \langle f_i:i < \mu \rangle \text{ such that}$
\nl 
$f_i \in G_\omega$ and $\mu \le \lambda \Rightarrow \pi_{n,\omega}(f_i) =
\pi_{n,\omega}(f_0)$ for $i < \mu$ 
and $i < j < \mu \and t \in I \Rightarrow \neg f_i {\Cal E}_t f_j\}$. \nl
We write nu$({\Cal A}) =^+ \mu$ to mean that moreover the 
supremum is obtained.  Let nu$^+({\Cal A})$ is the first $\mu$ such that
there is no $\langle f_i:i < \mu \rangle$ as above (so nu$({\Cal A}) \le
\text{nu}^+({\Cal A}) \le \text{ nu}({\Cal A})^+$ and nu$({\Cal A}) <
\text{ nu}^+({\Cal A})$ implies nu$({\Cal A})$ is a limit cardinal and the
supremum not obtained). \nl
4) We say ${\Cal A}$ is an explicit $(\bar \lambda,\bar J)$-system if:
${\Cal A} = (\bar \lambda,\bold{\bar J},\bar G,\bar H,\bar \pi,\bar \sigma)$
and
\mr
\item "{$(\alpha)$}"  $\bar \lambda = \langle \lambda_n:n < \omega \rangle,
\bar{\bold J} = \langle \bold J_n:n < \omega \rangle$
\sn
\item "{$(\beta)$}"  $\lambda_n < \lambda_{n+1},\bold J_n \subseteq
\bold J_{n+1}$,
\sn
\item "{$(\gamma)$}"  letting $\lambda^{\Cal A} = \dsize \sum_{n < \omega}
\lambda_n,\bold I^{\Cal A} = \dbcu_{n < \omega} \bold J_n$ we have
sys$({\Cal A}) =: (\lambda,\bold I,\bar G,\bar H,\bar \pi,\bar \sigma)$ is a
$(\lambda,\bold I)$-system
\sn
\item "{$(\delta)$}"  $|\bold J_n| \le \lambda_n,|G_n| \le \lambda_m,
|H^t_n| < \lambda$ and $|H^n_t| \le |H^{n+1}_t|$.
\ermn
5) We add in (4), full if
\mr
\item "{$(\varepsilon)$}"  $|H^t_n| \le \lambda_n$.  
\ermn
Second assume cf$(\mu) = \aleph_0$, so let $\mu = \dsize \sum_{n < \omega}
\mu_n,\mu_n< \mu_{n+1}$, and without loss of generality 
$\lambda_n < \mu_n = \text{ cf}(\mu_n)$
and $\mu > \lambda \Rightarrow \mu_n > \lambda$.  If $\mu > \lambda$, for each
$n$ there is a witness $\langle f^n_\alpha:\alpha < \mu_n \rangle$ to
nu$^+({\Cal A}) > \mu_n$, so $f^n_\alpha \in G^{\Cal A}_\omega$ and as
$\mu_n > \lambda \ge |G^{\Cal A}_n|$, \wilog \, $\pi_{n,\omega}(f^n_\alpha) =
\pi_{n,\omega}(f^0_\alpha)$.  Now for some increasing $\eta \in {}^\omega
\omega$ we have $n < \omega \and \alpha < \mu_n \Rightarrow \pi_{n,\omega}
(f^n_\alpha) \in G_n$ and we can continue as above. \nl
Lastly, assume $\mu = \lambda$, then rereading Definition \scite{5.1}(3), we
can choose $\mu_n = \lambda^+_n,\langle f^n_\alpha:\alpha < \mu_n \rangle$
as above, and continue as above.
\enddefinition
\bigskip

\proclaim{\stag{5.2} Claim}  1) For any strict $(\lambda,\bold I)$-system
${\Cal A}$ there is an explicit $(\bar \lambda,\bar{\bold J})$-system 
${\Cal B}$ such that sys$({\Cal B}) = {\Cal A}$ so

$$
\lambda = \dsize \sum_{n < \omega} \lambda_n,\bold I = \dbcu_{n < \omega}
\bold J_n, \text{ nu}({\Cal B}) = \text{nu}({\Cal A})
$$
\mn
(and if in one side the supremum is obtained, so in the other). \nl
2) For any $(\lambda,\bold I)$-system ${\Cal A}$ such that $\lambda >
2^{\aleph_0}$ and nu$^+({\Cal A}) > \mu$, cf$(\mu) \notin [\aleph_1,
2^{\aleph_0}]$ there is an explicit $(\bar \lambda,\bar{\bold J})$-system
${\Cal B}$ such that $\lambda^{\Cal A} = \dsize \sum_{n < \omega}
\lambda^{\Cal B}_n,\bold I^{\Cal A} = \dbcu_{n < \omega} \bold J^{\Cal B}_n$
and nu$^+({\Cal B}) > \mu$. \nl
3) In part (2) if $f:\text{Card } \cap \lambda \rightarrow$ Card is increasing
we can demand $\lambda_n \in \text{ Rang}(f)$, \nl
$f(\lambda_n) < \lambda_{n+1}$.
So if $\lambda$ is strong limit $> \aleph_0$, \ub{then} we can demand
$2^{\lambda^{\Cal B}_n} < \lambda^{\Cal B}_{n+1} = \text{cf}
(\lambda^{\Cal B}_{n+1})$.
\endproclaim
\bigskip

\demo{Proof}  1) Straight. \nl
2) Let $\bar \lambda = \langle \lambda_n:n < \omega \rangle$ be such that
$\lambda = \dsize \sum_{n < \omega} \lambda_n,2^{\aleph_0} < \lambda_n <
\lambda_{n+1}$, cf$(\lambda_n) = \lambda_n$.  Let $\langle G_{n,\ell}:\ell
< \omega \rangle$ be increasing, $G_{n,\ell}$ a subgroup of $G_n$ of
cardinality $\le \lambda_\ell$ and $G_n = \dbcu_{\ell < \omega} G_{n,\ell}$.
Let $\langle H^t_{n,\ell}:\ell < \omega \rangle$ be an increasing sequence of
subgroups of $H^t_n$ with union $H^t_n,|H^t_{n,\ell}| \le \lambda_\ell$.
Let $\langle \bold J_n:n < \omega \rangle$ be an increasing sequence of
subsets of $\bold I$ with union $\bold I$ such that $|\bold J_n| \le
\lambda_n$.
\sn
Without loss of generality $\pi_{m,n}$ maps $G_{n,\ell}$ into $G_{m,\ell}$
and $\pi^t_{m,n}$ maps $H^t_{n,\ell}$ into $H^t_{m,\ell}$ and
$\sigma^t_n$ maps $H^t_{n,\ell}$ into $G^t_{n,\ell}$ (why?  just chose
witness).
\sn
Now for every increasing $\eta \in {}^\omega \omega$ we let

$$
G^\eta_\omega = 
\{g \in G_\omega:\text{ for every } n < \omega \text{ we have }
\pi_{n,\omega}(g) \in G_{n,\eta(n)}\}.
$$
\mn
Clearly
\mr
\widestnumber\item{$(*)_1(\alpha)$}
\item "{$(*)_1(\alpha)$}"  $G^\eta_\omega$ 
is a subgroup of $G_{\omega,\eta}$ and
\sn
\item "{$(\beta)$}"  $\{G^\eta_\omega:
\eta \in {}^\omega \omega \text{ increasing}\}$ is directed
\sn
\item "{$(\gamma)$}"  $G_\omega = 
\cup\{G^\eta_\omega:\eta \in {}^\omega \omega \text{ (increasing)}\}$.
\ermn
First assume cf$(\mu) \ne \aleph_0$ so as 
cf$(\mu) > 2^{\aleph_0}$ for some $\eta \in {}^\omega \omega$, strictly
increasing, we have
\mr
\item "{$(*)_2$}"  nu$({\Cal A})$ is equal to 
sup$\{|X|:X \subseteq G_{\omega,\eta}$
and $t \in \bold I \and f \ne g \in X \Rightarrow f g^{-1} \notin
\sigma^t_\omega(H^t_\omega)\}$ \nl
(moreover, if one supremum is attained then so does the other).
\ermn
For each strictly increasing $\nu \in {}^\omega \omega$ let
$H^{(t,\nu)}_\omega$ be the subgroup $\{g \in H^t_\omega:\text{for every}$
\nl
$n < \omega \text{ we have } \sigma_{n,\omega}(g) \in H^{(t,\nu)}_\alpha\}$.
So let $\bold J'_n = \{(t,\nu):t \in \bold J \text{ and } \nu \in
{}^\omega \omega \text{ increasing}\}$.
\sn
We define $G_{n,\eta(n),\zeta}$, a subgroup of $G_{n,\eta(n)}$, decreasing
with $\zeta$ by induction on $\zeta$:
\mn
\ub{$\zeta=0$}:  $G_{n,\eta(n),\zeta} = G_{n,\eta(n)}$
\bn
\ub{$\zeta = \varepsilon + 1$}: $G_{n,\eta(n),\zeta} = \{x:x \in
G_{n,\eta(n),\varepsilon}$ and $x \in \text{ Rang}(\pi_{n,n+1} \restriction
G_{n+1,\eta(n+1),\varepsilon})$ \nl

$\qquad \qquad \qquad \qquad \qquad \quad$ and 
$n > 0 \Rightarrow \pi_{n-1,n}(x) \in
G_{n-1,\eta(n-1),\varepsilon}\}$
\bn
\ub{$\zeta$ limit}:  $G_{n,\eta(n),\zeta} = \dbca_{\varepsilon < \zeta}
G_{n,\eta(n),\varepsilon}$.
\bigskip

Let $G^\eta_n = \dbca_{\zeta < \lambda^+} G_{n,\eta(n),\zeta},
\pi^\eta_{m,n} = \pi_{m,n} \restriction G^\eta_n$.  Easily $\langle G^\eta_n,
\pi^\eta_{m,n}:
m \le n < \omega \rangle$ is directed with limit $G^\eta_\omega$
with $\pi^\eta_{n,\omega} = \pi_{n,\omega} \restriction G^\eta_\omega$.
\mn
Define $H^{(t,\nu)}_n,\pi^{(t,\nu)}_{m,m}$ parallely 
to $G^\eta_n,\pi^\eta_{m,n}$
but such that $\sigma^t_\alpha$ maps $H^{(t,\nu)}_\alpha$ into 
$G^\eta_\alpha$ (note: element of $H^{(t,\nu)}_\alpha$ not mapped to
$G^\eta_\alpha$ are irrelevant). \nl
Let $\sigma^{(t,\nu)}_\alpha:H^{(t,\nu)}_\alpha \rightarrow G^\eta_\omega$ be
$\sigma^t_\alpha \restriction H^{(t,\nu)}_\omega$ and
$\sigma^{(t,\sigma)}_n = \sigma^t_n \restriction H^{(t,\nu)}_n$.

We have defined actually $\bar B = (\bar \lambda^{\Cal B},\bold{\bar J}^
{\Cal B},\bar G,\bar H,\bar \pi^{\Cal B},\bar \sigma^{\Cal B})$ where \nl
$\bar \lambda^{\Cal B} = \langle \lambda_n:n < \omega \rangle,
\bold J^{\Cal B} = \langle \bold J_n \times {}^\omega \omega:n < \omega
\rangle,\bar G^{\Cal B} = \langle G^\eta_\alpha:\alpha \le \omega 
\rangle$, \nl
$\bar H^{\Cal B} = \left< \langle H^x_\alpha:\alpha \le \omega \rangle:
x \in \dbcu_n \bold J_n \right>$, \nl
$\bar \pi^{\Cal B} = \langle \pi^\eta_{\alpha,
\beta}:\alpha \le \beta \le \omega \rangle \char 94 \left< \langle
\pi^{(t,\nu)}_{\alpha,\beta}:\alpha \le \beta \le \omega \rangle:(t,\nu) \in
\dbcu_n \bold J_n \right>$ and \nl
$\bar \sigma^{\Cal B} = \left< \langle
\sigma^{(t,\nu)}_\alpha:\alpha \le \omega \rangle:(t,\nu) 
\in \dbcu_{n < \omega} \bold J_n \right>$.

We have almost finished.  Still $G^\eta_n$ may be of cardinality 
$> \lambda_n$ but note that for $k:\omega \rightarrow \omega$ non-decreasing
with limit $\omega,\langle G^\eta_n:n < \omega \rangle$ can be replaced by
$\langle G_{k(n)}:n < \omega \rangle$. \hfill$\square_{\scite{5.2}}$
\enddemo
\bn
For the rest of this section we adopt:
\demo{\stag{5.3A} Convention}  1) ${\Cal A}$ is an explicit $(\bar \lambda,
\bar{\bold J})$-system, so below rk$_t(g,f)$ should be written as
rk$_t(g,f,{\Cal A})$, etc. \nl
2) $\lambda = \dsize \sum_{n < \omega}
\lambda_n,\lambda_n = \lambda^{\Cal A}_n,\bold J_n = \bold J^{\Cal A}_n,
\bold I = \bold I^{\Cal A} = \dbcu_{n < \omega} \bold J_n,
G_\alpha = G^{\Cal A}_\alpha$, etc.
\nl
3) $k_t(n) = \text{ Max}\{m:m \le n,|H^t_m| \le \lambda_n\}$ so $k_t:\omega 
\rightarrow \omega$ is non-decreasing converging to $\infty$. \nl
For the reader's convenience we repeat \scite{5.5} - \scite{5.8} from $Gr$.
\enddemo
\bigskip

\definition{\stag{5.4} Definition}  1) For $g \in H^t_\alpha$ let lev$(g) 
= \alpha$ (\wilog \, this is well defined). \nl
2) For $\alpha \le \beta \le \omega,g \in H^t_\beta$ let $g \restriction
H^t_\alpha = \pi^t_{\alpha,\beta}(g)$ and we say $g \restriction H^t_\alpha$
is below $g$ and $g$ is above $g \restriction H^t_\alpha$ or extend
$g \restriction H^t_\alpha$. \nl
3) For $\alpha \le \beta \le \omega,f \in G_\beta$ let $f \restriction
G_\alpha = \pi_{\alpha,\beta}(f)$.
\sn
We will now describe the rank function used in the proof of the main
theorem.
\enddefinition
\bigskip

\definition{\stag{5.5} Definition}  1) For $g \in H^t_n,f \in G_\omega$ we
say that $(g,f)$ is a nice $t$-pair if $\sigma^t_n(g) = f \restriction G_n$.
\nl
2)  Define a ranking function rk$_t(g,f)$ for any nice $t$-pair.  
First by induction on the ordinal $\alpha$ (we can fix $f \in G_\omega$), 
we define when rk$_t(g,f) \ge \alpha$ simultaneously for all $n < \omega,
g \in H^t_n$
\mr
\item "{$(a)$}"  rk$_t(g,f) \ge 0$ \ub{iff} $(g,f)$ is a nice $t$-pair
\sn
\item "{$(b)$}"  rk$_t(g,f) \ge \delta$ for a limit ordinal $\delta$ \ub{iff}
for every $\beta < \delta$ we have rk$_t(g,f) \ge \beta$
\sn
\item "{$(c)$}"  rk$_t(g,f) \ge \beta + 1$ \ub{iff} 
$(g,f)$ is a nice $t$-pair,
and letting $n = \text{ lev}(g)$ there exists $g' \in H^t_{n+1}$ extending
$g$ such that rk$_t(g',f) \ge \beta$
\sn
\item "{$(d)$}"  rk$_t(g,f) \ge -1$.
\ermn
3) For $\alpha$ an ordinal or $-1$ (stipulating $-1 < \alpha < \infty$ for
any ordinal $\alpha$) we have rk$_t(g,f) = \alpha$ \ub{iff} 
rk$_t(g,f) \ge \alpha$ and it is false that rk$_t(g,f) \ge \alpha +1$. \nl
4) rk$_t(g,f) = \infty$ iff for every ordinal $\alpha$ we have
rk$_t(g,f) \ge \alpha$. 
\sn
The following two claims give the principal properties of rk$_t(g,f)$.
\enddefinition
\bigskip

\proclaim{\stag{5.6} Claim}  Let $(g,f)$ be a nice $t$-pair. \nl
1) The following statements are equivalent:
\mr
\item "{$(a)$}"  rk$_t(g,f) = \infty$
\sn
\item "{$(b)$}"  there exists $g' \in H^t_\omega$ extending $g$ such that
$\sigma^t_\omega(g') = f$.
\ermn
2)  If rk$_t(g,f) < \infty$, \ub{then} rk$_t(g,f) < \mu^+$ where
$\mu = \dsize \sum_{n < \omega} 2^{\lambda_n}$ (for $\lambda$ strong limit,
$\mu = \lambda$). \nl
3)  If $g'$ is a proper extension of $g$ and $(g',f)$ is also a nice $t$-pair
\ub{then}
\mr
\item "{$(\alpha)$}"  rk$_t(g',f) \le \text{ rk}_t(g,f)$ and
\sn
\item "{$(\beta)$}"  if $0 \le \text{ rk}_t(g,f) < \infty$ \ub{then} 
the inequality is strict.
\endroster
\endproclaim
\bigskip

\demo{Proof}  \nl
1) \ub{Statement $(a) \Rightarrow (b)$}. \nl
Let $n$ be the value
such that $g \in H^t_n$.  If we will be able to define $g_k \in H^t_k$ for
$k < \omega,k \ge n$ such that  
\mr
\widestnumber\item{$(iii)$}
\item "{$(i)$}"  $g_n = g$
\sn
\item "{$(ii)$}"  $g_k$ is below $g_{k+1}$ that is $\pi^t_{k,k+1}
(g_{k+1}) = g_k$ and
\sn
\item "{$(iii)$}"  rk$_t(g_k,f) = \infty$,
\ermn
\ub{then} clearly we will be done
since $g' =: \underset k {}\to {\text{lim}} \, g_k$ is as required.  The 
definition is by induction on $k \ge n$. \nl
For $k = n$ let $g_0 = g$. \nl
For $k \ge n$, suppose $g_k$ is defined.  By $(iii)$ we have rk$_t(g_k,f) =
\infty$, hence for every ordinal $\alpha$, rk$_t(g,f) > \alpha$ hence there
is $g^\alpha \in H^t_{k+1}$ extending $g$ such that rk$_t(g^\alpha,f) \ge
\alpha$.  Hence there exists $g^* \in H^t_{k+1}$ extending $g_k$ such that
$\{\alpha:g^\alpha = g^*\}$ is unbounded hence
rk$_t(g^*,f) = \infty$, and let $g_{k+1} =: g^*$.
\bn
\ub{Statement $(b) \Rightarrow (a)$}.  \nl
Since $g$ is below $g'$, it is enough to prove by induction on $\alpha$ that
for every $k \ge n$ when $g_k =: g' \restriction H^t_k$ we have that
rk$_t(g,f) \ge \alpha$.

For $\alpha = 0$, since $\sigma^t_\omega(g') = f \restriction G_n$ clearly for
every $k$ we have $\sigma^t_k(g_k) = f \restriction G_k$ so $(g_k,f)$ is a
nice $t$-pair.

For limit $\alpha$, by the induction hypothesis for every $\beta < \alpha$
and every $k$ we have rk$_t(g_k,f) \ge \beta$, hence by Definition 
\scite{5.5}(2)(b), rk$_t(g_k,f) \ge \alpha$.

For $\alpha = \beta +1$, by the induction hypothesis for every $k \ge n$ 
we have rk$_t(g_k,f) \ge \beta$.  Let $k_0 \ge n$ be given.  Since
$g_{k_0}$ is below $g_{k_0+1}$ and rk$_t(g_{k_0+1},f) \ge \beta$,
Definition \scite{5.5}(2)(c) implies that rk$_t(g_{k_0},f) \ge \beta +1$; i.e.
for every $k \ge n$ we have rk$_t(g_k,f) \ge \alpha$.  So we are done. \nl
2) Let $g \in H^t_n$ and $f \in G_\omega$ be given.  It is enough to prove
that if rk$_t(g,f) \ge \lambda^+$ then rk$_t(g,f) = \infty$.  Using part (1)
it is enough to find $g' \in H^t_\omega$ such that $g$ is below $g'$ and
$f = \sigma^t_\omega(g')$.

We define by induction on $k < \omega,g_k \in H^t_{n+k}$ such that
$g_k$ is below $g_{k+1}$ and rk$_t(g_k,f) \ge \mu^+$.  For $k=0$ let
$g_k = g$.  For $k+1$, for every $\alpha < \mu^+$, as rk$_t(g_k,f) >
\alpha$ by \scite{5.5}(2)(c) there is $g_{k,\alpha} \in G_{n+k+1}$ extending
$g_k$ such that rk$_t(g_{k,\alpha},f) \ge \alpha$.  But the number of possible
$g_{k,\alpha}$ is $\le |H^t_{n+k+1}| \le 2^{\lambda_n} +k+1 < \mu^+$
hence there are a function $g$ and a set $S \subseteq \mu^+$ of 
cardinality $\lambda^+$ such that $\alpha \in S \Rightarrow g_{k,\alpha}
=g$.  Then take $g_{k+1} = g$. \nl
3) Immediate from the definition.  \hfill$\square_{\scite{5.6}}$
\enddemo
\bigskip

\proclaim{\stag{5.7} Lemma}  1) Let $(g,f)$ be a nice $t$-pair.  \ub{Then}
we have rk$(g,f) \le \text{ rk}(g^{-1},f^{-1})$. \nl
2) For every nice $t$-pair $(g,f)$ we have rk$(g,f) = \text{ rk}(g^{-1},
f^{-1})$.
\endproclaim
\bigskip

\demo{Proof}  1) By induction on $\alpha$ prove that rk$(g,f) \ge \alpha
\Rightarrow \text{ rk}(g^{-1},f^{-1}) \ge \alpha$ (see more details in
Lemma \scite{5.8}). \nl
2) Apply part (1) twice. \hfill$\square_{\scite{5.7}}$
\enddemo
\bigskip

\proclaim{\stag{5.8} Lemma}  1) Let $n < \omega$ be fixed, and let $(g_1,f_1),
(g_2,f_2)$ be nice $t$-pairs with $g_\ell \in H^t_n(\ell = 1,2)$.  \ub{Then} 
$(g_1g_2,f_1f_2)$ is a nice pair and rk$_t(g_1g_2,f_1f_2) \ge 
\text{Min}\{\text{rk}_t(g_\ell,f_\ell):\ell =1,2\}$. \nl
2) Let $n,(f_1,g_1)$ and $(f_2,g_2)$ be as above.  If rk$_t(g_1,f_1) \ne
\text{ rk}_t(g_2,f_2)$, \ub{then} \nl
rk$_t(g_1g_2,f_1f_2) = \text{ Min}
\{\text{rk}_t(g_\ell,f_\ell):\ell =1,2\}$.
\endproclaim
\bigskip

\demo{Proof}  1) It is easy to show that the pair is $t$-nice.  
We show by induction
on $\alpha$ simultaneously for all $n < \omega$ and every $g_1,g_2 \in H^t_n$
that Min$\{\text{rk}(g_\ell,f_\ell):\ell =1,2\} \ge \alpha$ implies that
rk$(g_1g_2,f_1f_2) \ge \alpha$.

When $\alpha = 0$ or $\alpha$ is a limit ordinal this should be clear.  
Suppose $\alpha = \beta + 1$ and that rk$(g_\ell,f_\ell) \ge \beta +1$; by the
definition of rank for $\ell =1,2$ there exists $g'_\ell \in H^t_{n+1}$ 
extending $g_\ell$ such that $(g'_\ell,f_\ell)$ is a nice pair and 
rk$_t(g'_\ell,f_\ell) \ge \beta$.  By the induction assumption 
rk$_t(g'_1g'_2,f_1f_2) \ge \beta$.
Hence $g'_1g'_2$ is as required in the definition of rk$_t(g_1g_2,f_1f_2) 
\ge \beta + 1$. \nl
2) Suppose \, \wilog \, that rk$(g_1,f_1) < \text{ rk}(g_2,f_2)$, let
$\alpha_1 = \text{ rk}(g_1,f_1)$ and let $\alpha_2 = \text{ rk}_t(g_2,f_2)$.
By part (1), rk$_t(g_1g_2,f_1f_2) \ge \alpha_1$, by Proposition \scite{5.7},
rk$_t(g^{-1}_2,f^{-1}_2) = \alpha_2 > \alpha_1$.  So we have

$$
\align
\alpha_1 &= \text{ rk}_t(g_1,f_1) = \text{ rk}_t(g_1g_2g^{-1}_2,f_1f_2
f^{-1}_2) \\
  &\ge \text{ Min}\{\text{rk}_t(g_1g_2,f_1f_2),\text{rk}_t(g^{-1}_2,
f^{-1}_2)\} \\
  &= \text{ rk}_t(g_1g_2,f_1f_2) \ge \alpha_1.
\endalign
$$
\mn
Hence the conclusion follows.  \hfill$\square_{\scite{5.8}}$
\enddemo
\bigskip

\proclaim{\stag{5.9} Theorem}  Assume
\mr
\item "{$(a)$}"  $\lambda$ is strong limit $\lambda > \text{ cf}(\lambda)
= \aleph_0$
\sn
\item "{$(b)$}"  nu$({\Cal A}) \ge \lambda$.
\ermn
\ub{Then} nu$({\Cal A}) =^+ 2^\lambda$.
\endproclaim
\bn
The proof is broken into parts. \nl
\ub{\stag{5.10} Fact}:  Choose by induction on $n,\langle f_{n,i}:i < \lambda_n
\rangle$ such that
\mr
\item "{$(\alpha)$}"  $f_{n,i} \in G_\omega$ and $f_{n,i} \restriction
G_{n+1} = e_{G_{n+1}}$
\sn
\item "{$(\beta)$}"  $i < j < \lambda_n \and t \in \bold I \Rightarrow
\neg f_{n,i} {\Cal E}_t f_{n,j}$
\sn
\item "{$(\gamma)$}"  rk$_t(e_{H^t_k},f_{n,i}f^{-1}_{n,j}) < \infty$ for any
$t \in \bold J_n$ and $i \ne j < \lambda_n$ and $k = k_t$ (follows from
clause $(\beta)$)
\sn
\item "{$(\delta)$}"  if $f^*$ belongs to the subgroup of $G_\omega$
generated by the $\{f_{m,j}:m < n,j < \lambda_m\}$ and $t \in \bold J_n,
g \in \dbcu_{m \le k_t(n)} H^t_{k_t(n)}$, \ub{then} for every 
$i_0 < i_1 < i_2 < i_3 < \lambda_n$ each of the following statements 
have the same truth value, i.e. the truth value does not depend on
$(i_0,i_1,i_2,i_3)$)
{\roster
\itemitem{ $(i)$ }  rk$_t(g,f_{n,i_1}f^{-1}_{n,i_0}f^* f_{n,i_2} f^{-1}
_{n,i_3}) < \infty$
\sn
\itemitem{ $(ii)$ }  rk$_t(g,f_{n,i_3}f^{-1}_{n,i_2} f^* f_{n,i_0}
f^{-1}_{n,i_1}) < \infty$
\sn
\itemitem{ $(iii)$ }  rk$_t(e_{H^t_{k_t(n)}},f_{n,i_1}f^{-1}_{n,i_0}) <
\text{ rk}_t(g,f^*)$
\sn
\itemitem{ $(iv)$ }  rk$_t(e_{H^t_{k_t(n)}},f_{n,i_1}f^{-1}_{n,i_0}) >
\text{ rk}_t(g,f^*)$
\sn
\itemitem{ $(v)$ }  rk$_t(g,f^*) < \text{ rk}_t(g,f_{n,i_0}f^{-1}_{n,i_1}
f^* f_{n,i_2} f^{-1}_{n,i_3})$
\sn
\itemitem{ $(vi)$ }  rk$_t(g,f^*) < \text{ rk}_t(g,f_{n,i_2}f^{-1}_{n,i_3}
f^* f_{n,i_0} f^{-1}_{n,i_1})$
\endroster}
\item "{$(\varepsilon)$}"  for each $t \in \bold J_n$ one of the following
occurs:
{\roster
\itemitem{ $(a)$ }  for $i_0 < i_1 \le i_2 < i_3 < \lambda_n$ we have \nl
rk$_t(e_{H^t_{k_t(n)}},f_{n,i_0}f^{-1}_{n,i_1}) < \text{ rk}(e_{H^t_{k_t(n)}},
f_{n,i_2} f^{-1}_{n,i_3})$
\sn
\itemitem{ $(b)$ }  for some $\gamma^n_t$ for every $i < j < \lambda_n$
we have \nl
$\gamma^n_t = \text{ rk}_t(e_{H^t_{k_t(n)}},f_{n,i}f^{-1}_{n,j})$.
\endroster}
\endroster
\bigskip

\demo{Proof}  We can satisfy clauses $(\alpha),(\beta)$ by the definitions
and clause $(\gamma)$ follows.  Now clause $(\delta)$ is straight by 
Erd\"os Rado Theorem applied to a higher $n$. \nl
For clause $(\varepsilon)$ notice the transitivity of the order and of
equality.  \hfill$\square_{\scite{5.10}}$
\enddemo
\bigskip

\demo{\stag{5.14} Notation}  For $\alpha \le \omega$ let $T_\alpha = 
\times_{k < \alpha}\lambda_k,T =: \dbcu_{n < \omega} T_n$ (note: treeness
used).
\enddemo
\bigskip

\definition{\stag{5.17} Definition}  Now by 
induction on $n < \omega$, for every 
$\eta \in \dsize \prod_{m < n} \lambda_m$ we define $f_\eta \in G_\omega$
as follows:
\mr
\item "{{}}"  \ub{for $n=0$}:  $f_\eta = f_{<>} = e_{G_\omega}$
\sn
\item "{{}}"  \ub{for $n = m+1$}:  $f_\eta = f_{m,2 \eta(m)+1}
f^{-1}_{m,2 \eta(m)} f_{\eta \restriction m}$.
\endroster
\enddefinition
\bigskip

\demo{\stag{5.18} Fact}   For $\eta\in T_\omega$ and $m \le n < \omega$ we
have

$$
f_{\eta \restriction n} \restriction G_{m+1} = f_{\eta \restriction m}
\restriction G_{m+1}.
$$
\enddemo
\bigskip

\demo{Proof}   As $\pi_{n,\omega}$ is a homomorphism it is enough to prove
$(f_{\eta \restriction n}(f_{\eta \restriction m})^{-1}) \restriction
G_{n+1} = e_{G_{n+1}}$, hence it is enough to prove $n \le k < \omega 
\Rightarrow (f_{\eta \restriction k} f^{-1}_{\eta \restriction (k+1)}) 
\restriction G_{n+1} = e_{G_{n+1}}$ which follows from 
$k < \omega \Rightarrow f_{\eta \restriction k} 
f^{-1}_{\eta \restriction (k+1)} \restriction G_{k+1} = e_{G_{k+1}}$, 
which follows from $f_{k,\eta(\zeta)} 
\restriction G_{k+1} = e_{G_{k+1}}$ which 
holds by clause $(\alpha)$ above.  \hfill$\square_{\scite{5.18}}$
\enddemo
\bigskip

\definition{\stag{5.19} Definition}   For 
$\eta \in T_\omega$ we have $f_\eta \in G_\omega$ is well defined as
the inverse limit of 
$\langle f_{\eta \restriction n} \restriction G_n:n < \omega \rangle$, so 
$n < \omega \rightarrow f_\eta \restriction G_n =
f_{\eta \restriction n}$. Follows by \scite{5.18} and 
$G^\omega$ being an inverse limit. 
\enddefinition
\bigskip

\demo{\stag{5.20} Proposition}  Let $\eta,\nu \in T_\omega$ be such that
$(\forall^\infty n)(\eta(n) \ne \nu(n)),\eta(n) > 0,\nu(n) > 0$.  If
$t \in \bold I$, \ub{then} $f_\eta f^{-1}_\nu \notin \sigma^t_\omega
(H^t_\omega)$.
\enddemo
\bigskip

\demo{Proof}  Suppose toward contradiction that for some $g \in H^t_\omega$
we have $\sigma^t_\omega(g) = f_\eta f^{-1}_\nu$. \nl
Let $k < \omega$ be large enough such that $t \in \bold J_k,(\forall \ell)
[k \le \ell < \omega \rightarrow \eta(\ell) \ne \nu(\ell)]$.  Let
$\xi^\ell = \text{ rk}_t(g \restriction H^t_{k_t(\ell)},f_{\eta \restriction
(\ell +1)}f^{-1}_{\nu \restriction (\ell +1)})$ and $\zeta^\ell =
\text{ rk}_t(g \restriction H^t_{k_t(\ell +1)},f_{\eta \restriction (\ell +1)}
f^{-1}_{\nu \restriction (\ell +1)})$ (the difference between the two is 
the use of $k_t(\ell)$ via $k_t(\ell + 1))$.  Clearly
\mr
\item "{$(*)_1$}"  $f_{\eta \restriction (\ell +1)}f^{-1}_{\nu \restriction
(\ell +1)} = (f_{\ell,2 \eta(\ell)+1} f^{-1}_{\ell,2 \eta(\ell)})
(f_{\eta \restriction \ell} f^{-1}_{\nu \restriction \ell})
f_{\ell,2 \nu(\ell)}f^{-1}_{\ell,2 \nu(\ell)+1}$
\ermn
[Why?  Algebraic computations.]  Next we claim that
\mr
\item "{$(*)_2$}"  $\xi^\ell < \infty$ for $\ell \ge k$ ($\ell < \omega$).
\ermn
Why? \nl
\ub{Case 1}:  $\eta(\ell) < \nu(\ell)$.

Assume toward contradiction $\xi^\ell = \infty$, but 
by clause $(\gamma)$ above
\nl
rk$_t(e_{H^t_{k_t(\ell)}},f_{\ell,2 \eta(\ell)+2} 
f^{-1}_{\ell,2 \eta(\ell)+1}) < \infty = \xi^\ell$, hence by 
\scite{5.8}(3).

$$
\align
\text{rk}_t(g \restriction H^t_{k_t(\ell)},f_{\ell,2 \eta(\ell)+2} 
f^{-1}_{\ell,2 \eta(\ell)+1} f_{\eta \restriction (\ell+1)}
f^{-1}_{\nu \restriction (\ell +1)}) = \text{ Min}\{&\text{rk}_t
(e_{H^t_{k_t(\ell)}},f_{\ell,2(\eta(\ell)+2} f^{-1}_{\ell,2 \eta(\ell)+1}), \\
  &\text{rk}_t(g \restriction H^t_{k_t(\ell)},f_{\eta \restriction(\ell +1)}
f^{-1}_{\nu \restriction (\ell +1)})\} = \\
  &\text{rk}_t(e_{H^t_{k_t(\ell)}},f_{\ell,2 \eta(\ell)+2}
f^{-1}_{\ell,2 \eta(\ell)+1}) < \infty.
\endalign
$$
\mn
Now (by the choice of $f_{\eta \restriction (\ell +1)},
f_{\nu \restriction (\ell +1)}$, algebraic computation and 
the previous inequality) we have

$$
\align
\infty > \text{ rk}_t(g \restriction H^t_{k_t(\ell)},&f_{\ell,2 \eta(\ell)+2}
f^{-1}_{\ell,2 \eta(\ell)+1} f_{\eta \restriction (\ell +1)}
f^{-1}_{\nu \restriction (\ell+1)}) = \\
  &\text{rk}_t(g \restriction H^t_\ell,(f_{\ell,2 \eta(\ell)+2}
f^{-1}_{\ell,2 \eta(\ell)})(f_{\eta \restriction \ell}
f^{-1}_{\nu \restriction \ell})(f_{\ell,2\nu(\ell)}
f^{-1}_{\ell,2 \nu(\ell)+1})).
\endalign
$$
\mn
This and the assumption $\xi_\ell = \infty$ gives a contradiction to
$(\delta)(i)$ of \scite{5.10} (for $(i_0,i_1,i_2,i_3)$ being
$(2 \eta(\ell),2 \eta(\ell) +2,2 \nu(\ell),2 \nu(\ell)+1$) and being
$(2 \eta(\ell),2 \eta(\ell) +1,2 \nu(\ell),2 \nu(\ell)+1)$.
\bn
\ub{Case 2}:  $\nu(\ell) > \eta(\ell)$.

Similar using $(\delta)(ii)$ of \scite{5.10} instead of $(\delta)(i)$ of
\scite{5.10} (using $\eta(\ell) > 0$). \nl
So we have proved $(*)_2$.
\mr
\item "{$(*)_3$}"  $\xi^{\ell +1} \le \zeta^\ell$ for $\ell > k$. 
\ermn
Why?  Assume toward contradiction $\xi^{\ell +1} > \zeta^\ell$.  \nl
Let
$f^* = f_{\eta \restriction (\ell +1)} f^{-1}_{\nu \restriction (\ell +1)}$,
so $\zeta^\ell = \text{ rk}_t(g \restriction H^t_{k_t(\ell +1)},f^*)$ and 
using the choie of $\xi^{\ell +1}$ and $(*)_1$ we have
$\xi^{\ell +1} = \text{ rk}_t(g \restriction H^t_{k_t(\ell+1)},
f_{(\ell+1),2 \eta(\ell +1)+1} f^{-1}_{\ell +1,2 \eta(\ell +1)} f^*
f_{\ell +1,2 \nu(\ell +1)}$ \nl
$f^{-1}_{\ell +1,2 \nu(\ell +1)+1})$.
\mn

If $\zeta^\ell < \text{ rk}_t(e_{H^t_{k_t(\ell +1)}},f_{\ell +1,2 \eta
(\ell +1)+1} f^{-1}_{\ell+1,2 \eta(\ell +1)})$ then by \scite{5.10}$(\delta)
(iii)$ also \nl
$\zeta^\ell < \text{ rk}_t(e_{H^t_{k_t(\ell +1)}},
f_{\ell +1,2 \nu(\ell +1)+1} f^{-1}_{\ell+1,2 \nu(\ell +1)})$ hence using
twice \scite{5.8}(2) we have first $\zeta^\ell = \text{ rk}_t(g \restriction
H^t_{k_t(\ell +1)},f_{\ell +1,2 \eta(\ell +1)+1} 
f^{-1}_{\ell+1,2 \eta(\ell +1)}f^*)$ and second
(using also \scite{5.7}(2))
$\zeta^\ell = \text{ rk}_t(g \restriction H^t_{k_t(\ell +1)},
f_{\ell +1,2 \eta(\ell +1)+1} f^{-1}_{\ell+1,2 \eta(\ell +1)} f^*
f_{\ell +1,2 \eta(\ell +1)} f^{-1}_{\ell +1,2 \eta(\ell +1)+1})$,
so by the second statement in the previous paragraph we get $\zeta_\ell
= \xi^{\ell +1}$ contradicting our temporary assumption toward contradiction.

Also if rk$_t(e_{H^t_{k_t(\ell +1)}},f_{\ell +1,2 \eta(\ell +1)+1} 
f^{-1}_{\ell+1,2 \eta(\ell +1)}) \ne \text{ rk}_t(e_{H^t_{k_t(\ell +1)}},
f_{\ell +1,2 \nu(\ell +1)+1} f^{-1}_{\ell+1,2 \nu(\ell +1)}$ 
then by \scite{5.10}$(\delta)(iii) + (iv)$
also $\zeta^\ell$ is not equal to those two ordinals so similarly to the
previous sentence, \scite{5.8}(2) gives 
$\xi^{\ell +1} = \text{ Min}\{\text{rk}_t(e_{H^t_{k_t(\ell +1)}},
f_{\ell +1,2 \eta(\ell +1)+1} 
f^{-1}_{\ell+1,2 \eta(\ell +1)})$, \nl
$\text{rk}_t(g \restriction H^t_{k_t(\ell +1)},f^*),
\text{rk}_t(e_{H^t_{k_t(\ell +1)}},f_{\ell +1,2 \nu(\ell +1)+1} 
f^{-1}_{\ell+1,2 \nu(\ell +1)})\}$ which is $\le \zeta^\ell$ so
$\xi^{\ell +1} = \zeta^\ell$, contradicting our assumption toward 
contradiction.
\mn
Together the case left is (remember \scite{5.7})
\mr
\item "{$\boxtimes$}"  $\zeta^\ell = 
\text{ rk}_t(g \restriction H^t_{k_t(\ell +1)},f^*) \ge
\text{ rk}_t(e_{H^t_{k_t(\ell +1)}},f_{\ell +1,2 \eta(\ell +1)+1} 
f^{-1}_{\ell+1,2 \eta(\ell +1)}) = \text{ rk}_t(e_{H^t_{k_t(\ell +1)}},
f_{\ell +1,2 \nu(\ell +1)+1} f^{-1}_{\ell+1,2 \nu(\ell +1)})$. 
\ermn
So in \scite{5.10}$(\varepsilon)$, for $n = \ell +1$, case (b) holds. \nl
As, toward contradiction we are assuming $\xi^{\ell +1} > \zeta^\ell$ in 
the proof of $(*)_3$ we get, by \scite{5.8} e.g. if 
$\eta(\ell +1) > \nu(\ell +1)$ that \nl
$\text{rk}_t(g \restriction H^t_{k_t(\ell +1)},f_{\ell +1,2 \eta(\ell +1)+2} 
f^{-1}_{\ell+1,2 \eta(\ell +1)} f_{\eta \restriction (\ell +1)}
f^{-1}_{\nu \restriction (\ell +1)} f_{\ell +1,2 \nu(\ell +1)}
f^{-1}_{\ell +1,2 \nu(\ell +1)+1}) =
\text{rk}_t(e_{H^t_{k_t(\ell +1)}},f_{\ell +1,2 \eta(\ell +1)+2} 
f^{-1}_{\ell+1,2 \eta(\ell +1)})$ but by (b) of \scite{5.10}$(\varepsilon)$
proved above the later is $\le \zeta^\ell < \xi^{\ell +1}
= \text{ rk}_t(g \restriction H^t_{k_t(\ell +1)},
f_{\ell +1,2 \eta(\ell +1)+1}f^{-1}_{\ell+1,2 \eta(\ell)}
f_{\eta \restriction (\ell +1)} f^{-1}_{\nu \restriction (\ell+1)}
f^{-1}_{\ell +1,2 \nu(\ell +1)})$ contradiction to \scite{5.10}$(\delta)(v)$.
If $\eta(\ell +1) < \nu(\ell +1)$ we use similarly $f_{\ell +1,2 \nu
(\ell +1)+2}f^{-1}_{\ell +1,2 \nu(\ell +1)}$.
So $(*)_3$ holds.
\mr
\item "{$(*)_4$}"  $\zeta^\ell \le \xi^\ell$ \nl
[Why?  Look at their definitions, as $g \restriction H^t_{k_t(\ell +1)}$ is
above $g \restriction H^t_{k_t(\ell)}$.  Now if $k_t(\ell),k_t(\ell +1)$
are equal trivial otherwise use \scite{5.6}(3).]
\sn
\item "{$(*)_5$}"  if $k_t(\ell +1) > k_t(\ell)$ then $\zeta^\ell < \xi^\ell$
(so $\xi^\ell > 0$) \nl
[Why?  Like $(*)_4$.]
\sn
\item "{$(*)_6$}"  $\xi^\ell \ge \xi^{\ell +1}$ and if $k_t(\ell +1) >
k_t(\ell)$ then $\xi^\ell > \xi^{\ell +1}$ \nl
[Why?  By $(*)_3 + (*)_4$ the first phrase, and $(*)_3 + (*)_5$ for the
second phrase.]
\ermn
So $\langle \xi^\ell:\ell \in [k,\omega) \rangle$ is non-increasing, and
not eventually constant, contradiction.  \nl
${{}}$ \hfill$\square_{\scite{5.20}}$
\enddemo
\bigskip

\demo{Proof of \scite{5.9}}  Obvious as we can find $T' \subseteq T$, a
subtree with $\lambda^{\aleph_0}$ $\omega$-branches and $\eta \ne \nu \in
\text{ lim}(T') \Rightarrow (\forall^\infty \ell)\eta(\ell) \ne \nu(\ell)$
and $\eta \in \lim(T') \and n < \omega \Rightarrow \eta(n) > 0$. \nl
Now $\langle f_\eta:\eta \in \lim(T') \rangle$ is as required.
\enddemo
\bn
\ub{\stag{5.21} Conclusion}:  If ${\Cal A}$ is a $(\lambda,\bold I)$-system,
and $\lambda$ is a strong limit of cofinality $\aleph_0$ and nu$({\Cal A})
\ge \lambda$, \ub{then} nu$({\Cal A}) =^+ 2^\lambda$.
\bigskip

\demo{Proof}  If $\lambda = \aleph_0$ by \scite{1.1}, if $\lambda > \aleph_0$
then $\lambda > 2^{\aleph_0}$.  We apply \scite{5.2}(2) to ${\Cal A}$ and
$\mu = \lambda$ (so cf$(\mu) = \aleph_0$) and get ${\Cal B}$, for which by
\scite{5.8} we have nu$({\Cal B}) = 2^\lambda$ hence by the choice of
${\Cal B}$ also nu$({\Cal B}) =^+ 2^\lambda$.
\hfill$\square_{\scite{5.21}}$
\enddemo
\newpage

\shlhetal

\newpage
    
REFERENCES.  
\bibliographystyle{lit-plain}
\bibliography{lista,listb,listx,listf,liste}

\def\germ{\frak} \def\scr{\cal}
  \ifx\documentclass\undefinedcs\def\rm{\fam0\tenrm}\fi
  \def\defaultdefine#1#2{\expandafter\ifx\csname#1\endcsname\relax
  \expandafter\def\csname#1\endcsname{#2}\fi} \defaultdefine{Bbb}{\bf}
  \defaultdefine{frak}{\bf} \defaultdefine{mathbb}{\bf}
  \defaultdefine{beth}{BETH} \def\bbfI{{\Bbb I}} \def\mbox{\hbox}
  \def\text{\hbox} \def\om{\omega} \def\Cal#1{{\bf #1}} \def\pcf{pcf}
  \defaultdefine{cf}{cf} \def\restriction{{|}} \def\club{CLUB} \def\w{\omega}
  \def\exist{\exists} \def\se{{\germ se}} \def\bb{{\bf b}}
  \def\equivalence{\equiv} \def\cite#1{[#1]} \def\implies{\Rightarrow}
\begin{thebibliography}{GrSh 302a}
\makeatletter \renewcommand{\@biblabel}[1]{[#1]} \makeatother

\bibitem[EkSh 505]{EkSh:505}Paul~C. Eklof and Saharon Shelah.
\newblock {A Combinatorial Principle Equivalent to the Existence of Non-free
  Whitehead Groups}.
\newblock In {\em Abelian group theory and related topics}, volume 171 of {\em
  {Contemporary Mathematics}}, pages 79--98. {American Mathematical Society,
  Providence, RI}, 1994.
\newblock edited by R. Goebel, P. Hill and W. Liebert, Oberwolfach proceedings.

\bibitem[GiSh 597]{GiSh:597}Moti Gitik and Saharon Shelah.
\newblock {On densities of box products}.
\newblock {\em {Topology and its Applications}}, {\bf accepted}.

\bibitem[GrSh 302a]{GrSh:302a}Rami Grossberg and Saharon Shelah.
\newblock {On cardinalities in quotients of inverse limits of groups}.
\newblock {\em {Mathematica Japonica}}, {\bf 47}(2), 1998.

\bibitem[MRSh 314]{MRSh:314}Alan~H. Mekler, Andrzej Ros{\l}anowski, and Saharon
  Shelah.
\newblock {On the $p$-rank of Ext}.
\newblock {\em {Israel Journal of Mathematics}}, {\bf submitted}.

\bibitem[MkSh 418]{MkSh:418}Alan~H. Mekler and Saharon Shelah.
\newblock {Every coseparable group may be free}.
\newblock {\em {Israel Journal of Mathematics}}, {\bf 81}:161--178, 1993.

\bibitem[Sh 273]{Sh:273}Saharon Shelah.
\newblock {Can the fundamental (homotopy) group of a space be the rationals?}
\newblock {\em {Proceedings of the American Mathematical Society}}, {\bf
  103}:627--632, 1988.

\bibitem[Sh:f]{Sh:f}Saharon Shelah.
\newblock {\em {Proper and improper forcing}}.
\newblock {Perspectives in Mathematical Logic}. {Springer}, 1998.

\end{thebibliography}

\enddocument

\bye